\documentclass[12pt]{article}
\usepackage{amsmath}
\usepackage{amsthm}
\usepackage{amssymb}
\usepackage{tabularx}
\usepackage{multirow}
\usepackage{tabularx}

\title{Isomorphisms of Kac-Moody groups which preserve bounded subgroups}
\author{Pierre\,-Emmanuel {\sc Caprace}\footnotemark{} \ and Bernhard {\sc M\"uhlherr}}
\date{February 2005}
\newtheorem{prop}{Proposition}[section]
\newtheorem{theo}[prop]{Theorem}
\newtheorem*{theo_intro}{Theorem}
\newtheorem{corintro}{Corollary}

\newtheorem{cor}[prop]{Corollary}
\newtheorem{lem}[prop]{Lemma}

\theoremstyle{remark}
\newtheorem{remark}[prop]{Remark}

\newcommand{\eps}{\epsilon}
\newcommand{\cod}{\delta^{\ast}}
\newcommand{\K}{\mathbb{K}}

\newcommand{\TRD}{$(G, (U_{\alpha})_{\alpha \in \Phi})$}
\newcommand{\tb}{((\Delta_+, \delta_+),( \Delta_-,
\delta_-), \delta^{\ast})}
\newcommand{\la}{\langle}
\newcommand{\ra}{\rangle}
\DeclareMathOperator{\proj}{proj} \DeclareMathOperator{\op}{op}
\DeclareMathOperator{\Stab}{Stab} \DeclareMathOperator{\Fix}{Fix}
\DeclareMathOperator{\Res}{Res}
\begin{document}
\maketitle
\begin{center}
{\small D\'epartement de Math\'ematiques, \\Universit\'e libre de Bruxelles, CP216\\
Bd du Triomphe, B-1050 Bruxelles, Belgium\\e-mail : {\it
pcaprace@ulb.ac.be} and {\it bmuhlherr@ulb.ac.be}}\\
\vspace{1cm}
\begin{tabularx}{12cm}{X}
{\sc Abstract. } A subgroup of a Kac-Moody group is called
\textbf{bounded} if it is contained in the intersection of two
finite type parabolic subgroups of opposite signs. In this paper,
we study the isomorphisms between Kac-Moody groups over arbitrary
fields of cardinality at least 4, which preserve the set of
bounded subgroups. We show that such an isomorphism between two
such Kac-Moody groups induces an isomorphism between the
respective twin root data of these groups. As a consequence, we
obtain the solution of the isomorphism problem for Kac-Moody
groups over finite fields of cardinality at least 4.
\end{tabularx}
\end{center}
\makeatletter
\renewcommand{\@makefnmark}{\mbox{$^{*}$}}%
\renewcommand{\@makefntext}[1]{\noindent\makebox[0.4 em][r]{\@makefnmark}#1}
\footnotetext{F.N.R.S. research fellow
} \makeatother
\bigskip
\bigskip
\section{Introduction}
Kac-Moody groups are infinite-dimensional generalizations of
Chevalley groups. It is known that each automorphism of a
Chevalley group (of irreducible type and over a perfect field) can
be written as a product of an inner, a diagonal, a graph and a
field automorphism (see Theorem 30 in \cite{St68}). In \cite{KP87}
it was conjectured that the same statement holds for Kac-Moody
groups over algebraically closed fields of characteristic 0 up to
the addition of a so called sign automorphism. In \cite{CM03b}
this conjecture is shown to be true for Kac-Moody groups over
algebraically closed fields of any characteristic. This is
achieved in loc. cit. by solving the isomorphism problem for those
groups. In this paper, we study the isomorphism problem for
Kac-Moody groups over arbitrary fields of cardinality at least
$4$. We restrict our attention to isomorphisms which preserve the
set of bounded subgroups. In this context, a subgroup of a
Kac-Moody group is called \textbf{bounded} if it is contained in
the intersection of two finite type parabolic subgroups of
opposite signs.

Throughout the paper we use Tits' definition for Kac-Moody groups
over fields \cite{Ti87}. This definition does not only provide the
abstract Kac-Moody group $G$ but also a canonical system
$(U_\alpha)_{\alpha \in \Phi}$ of root subgroups. The pair $(G,
(U_\alpha)_{\alpha \in \Phi})$ is an example of a so called twin
root datum. Twin root data have been introduced by Tits in order
to give suitable axioms for these pairs arising from his
definition of Kac-Moody groups.

Our main result says that if two Kac-Moody groups are isomorphic
via such an isomorphism, then the groups are of the same type and
defined over the same ground field. Here is a precise statement.

\begin{theo_intro} Let
$\mathcal{D}=(G, (U_\alpha)_{\alpha \in \Phi})$ and
$\mathcal{D'}=(G',(U'_{\alpha'})_{\alpha' \in \Phi'})$ be two twin
root data associated with two Kac-Moody groups of non-spherical
type over fields of cardinality at least 4. Let $\xi : G \to G'$
be a group isomorphism which maps bounded subgroups of $G$ to
bounded subgroups of $G'$. Then $\xi$ induces an isomorphism of
$\mathcal{D}$ to $\mathcal{D}'$.
\end{theo_intro}

We refer to Section \ref{sect:TRD:def} below for the definition of
an isomorphism between twin root data. Roughly speaking, it means
that $(\xi(U_\alpha))_{\alpha \in \Phi}$ is `nearly
$G'$-conjugate' to $(U'_{\alpha'})_{\alpha' \in \Phi'}$.

As it is the case in the paper \cite{CM03b}, the present work
makes crucial use of the theory of twin buildings. A group endowed
with a twin root datum is indeed naturally endowed with a strongly
transitive action on a twin building, and the combinatorial
properties of this action turn out to be the most appropriate tool
in studying Kac-Moody groups from our point of view. However, we
have tried to explain each crucial building-theoretic statement in
more classical terms, without making reference to the language of
buildings. We hope this will help the reader who is not familiar
to the theory of buildings to understand the main ideas of this
paper.

As a consequence of the theorem above, we obtain the following
result on automorphisms of Kac-Moody groups.

\begin{corintro}\label{cor 1}
Let $G$ be a Kac-Moody group over a field of cardinality at least
4. Let $\varphi$ be an automorphism of $G$ which preserves the set
of bounded subgroups. Then $\varphi$ splits as a product of an
inner, a diagonal, a graph, a field  and a sign automorphism.
\end{corintro}

There are mainly two motivations to consider isomorphisms which
preserve bounded subgroups.

The first motivation comes from the earlier work \cite{KW92} by
Kac and Wang. In this paper, automorphisms of Kac-Moody groups
over fields of characteristic 0 and associated with symmetrizable
Cartan matrices have been studied. One of the main results of loc.
cit. is that, given such a Kac-Moody group $G$ and its Kac-Moody
algebra $\mathfrak{g}$, then an automorphism of $G$ which
preserves the set of $\mathrm{Ad}_{\mathfrak{g}'}$-finite elements
splits as a product as in Corollary \ref{cor 1} above, where
$\mathfrak{g'}:=[\mathfrak{g, g}]$. We recall that an element $g
\in G$ is $\mathrm{Ad}_{\mathfrak{g}'}$-finite if and only if the
subgroup generated by $g$ is bounded (see \cite{KW92}, Theorem
2.10). Thus, Corollary~\ref{cor 1} can be seen as a weaker version
of Kac-Wang's result, which remains valid for Kac-Moody groups of
arbitrary type and over fields of arbitrary characteristic.

The second motivation is the fact that, in the case of a Kac-Moody
group over a finite field, a subgroup is bounded if and only if it
is finite (see Corollary~\ref{G_f <=> finite} below). Therefore,
all isomorphisms preserve bounded subgroups in this case.
Consequently, we obtain the following result.

\begin{corintro}\label{cor 2}
Let $\mathcal{D}=(G, (U_\alpha)_{\alpha \in \Phi})$ and
$\mathcal{D'}=(G',(U'_{\alpha'})_{\alpha' \in \Phi'})$ be two twin
root data associated with two Kac-Moody groups over finite fields
of cardinality at least 4 and let $\xi : G \to G'$ be an
isomorphism. Then $\xi$ induces an isomorphism of $\mathcal{D}$ to
$\mathcal{D'}$ unless $G$ and $G'$ are both finite. Moreover, any
automorphism of $G$ splits as a product of an inner, a diagonal, a
graph, a field  and a sign automorphism.
\end{corintro}

Kac-Moody groups over finite fields are finitely generated and
some subclasses of them are known to be finitely presented. In the
recent years these groups became important in geometric group
theory for several reasons (see \cite{Re04}). In this context, B.
R\'emy proved a factorization theorem for the automorphisms of
certain Kac-Moody groups (see \cite{Re02}, Theorem 3.1). Corollary
2 covers this result as a special case.

Let us also mention the existence of exotic constructions of
groups of Kac-Moody type, which are not Kac-Moody groups in the
strict sense but which are also endowed with a twin root datum.
For example, R\'emy and Ronan~\cite{RR} constructed examples of
groups of Kac-Moody type defined simultaneously over different
ground fields. It turns out that, provided the maximal tori are
locally large enough, our methods extend also to these exotic
cases, and the interested reader will have no difficulty to extend
our arguments to this slightly more general situation (see also
the introduction of \cite{CM03b} for other remarks and results
related to the isomorphism problem of exotic groups of Kac-Moody
type).

The paper is organized as follows. After a preliminary section
where definitions are recalled, notation is fixed and some
auxiliary results are proven, we discuss in Section \ref{Levi
sect} the Levi decomposition of the intersection of two parabolic
subgroups of opposite sign in a group endowed with a twin root
datum (a similar but slightly less general discussion had been
done in \cite{Re99}). The key result of this paper is contained in
Section \ref{max sub sect}, where we prove that the maximal
bounded subgroups coincide almost always with the Levi factors of
the maximal parabolic subgroups of finite type. In the next
section, we use this key result to state and prove a technical
version of our main result, which is valid for a larger class of
groups endowed with a twin root datum. Finally, the last two
sections are devoted to the proof of the main theorem above and
its corollaries.

\section{Preliminaries}\label{preliminaries}
The main references are \cite{Ab96}, \cite{Re99}, \cite{Ro89} and
\cite{Ti92}.

We start by fixing a general convention~:\\
\emph{The ordered pair $(W, S)$ is a Coxeter system and $\ell$
denotes the corresponding length function. For $J \subseteq S$ we
set $W_J := \langle J \rangle$ and we call $J$ \textbf{spherical}
whenever $W_J$ is finite. }

\subsection{Buildings}

\subsubsection{Definition}\label{sect:bldg:def}

A \textbf{building} of type $(W, S)$ is a set $\Delta$, whose
elements are called \textbf{chambers}, endowed with a map
$\delta~: \Delta \times \Delta \to W$ called the
\textbf{$W$-distance} satisfying the following axioms, where $x, y
\in \Delta$ and $w = \delta(x, y)$:
$$\begin{array}{ll}
\text{\textbf{(Bu1)}} & w = 1 \Leftrightarrow x = y;\\
\text{\textbf{(Bu2)}} & \text{if } z \in \Delta \text{ is such
that } \delta(y, z)= s \in S, \text{ then } \delta(x, z) \in \{w,
ws\};\\
& \text{if, furthermore, } \ell(ws) = \ell(w) + 1, \text{ then }
\delta(x, z) = ws;\\
\text{\textbf{(Bu3)}} & \text{if } s \in S, \text{ there exists }
z \in \Delta \text{ such that } \delta(y, z) = s \text{ and }
\delta(x, z)=ws.
\end{array}$$
For any two chambers $x, y \in \Delta$, the natural number
$\ell(\delta(x, y))$ is called the \textbf{numerical distance}
between $x$ and $y$.

An \textbf{isometry} between subsets of buildings of type $(W,S)$
is a bijection preserving the $W$-distance.

\subsubsection{Apartments}\label{sect:bldg:apt}

Given a Coxeter system $(W,S)$, let $\delta:W\times W\to W$ be
defined by $\delta:(x, y) \mapsto x^{-1}y$. In this way, we endow
$W$ with a canonical structure of a building of type $(W,S)$. This
building is denoted by $\mathcal{A}(W,S)$.

Any subset of a building $(\Delta, \delta)$ of type $(W,S)$ which
is isometric to the canonical building $\mathcal{A}(W,S)$ is
called an \textbf{apartment}. A fundamental property of buildings
is that any two chambers lie in a common apartment (see
\cite{Ro89}, Theorem 3.7).

\subsubsection{Panels, residues and
galleries}\label{sect:bldg:res}

Given $c \in \Delta$ and $s \in S$, then the set $\big\{x \in
\Delta | \delta(x, c) \in \{1, s\}\big\}$ is called an
\textbf{$s$-panel} of $\Delta$ or a \textbf{panel of type $s$}. A
\textbf{panel} is an $s$-panel for some $s \in S$. More generally,
for $c \in \Delta$ and $J \subseteq S$ the set
$$\mathrm{Res}_{J}(c):=\big\{x \in \Delta | \delta(x, c) \in W_J\}$$ is
called the \textbf{$J$-residue} of $\Delta$ which contains $c$.
Its \textbf{rank} is the cardinality of the set $J$; hence,
residues of rank~0 are just chambers and the residues of rank~1
are panels.

It is an important fact that a $J$-residue is itself a building of
type $(W_J, J)$ with the $W_J$-distance induced by $\delta$ (see
\cite{Ro89}, Theorem 3.5). Moreover, given a residue $R$ and an
apartment $\Sigma$ in a building $\Delta$, the intersection $R
\cap \Sigma$ is either empty or an apartment of $R$, and all
apartments of $R$ arise in this way. It is common and handy to say
that $R$ \emph{is contained in} $\Sigma$ whenever $R \cap \Sigma$
is nonempty.

A sequence of chambers such that two consecutive chambers are
\textbf{adjacent}, namely contained in a common panel, is called a
\textbf{gallery}. The gallery $\gamma = (x_0, x_1, \dots, x_n)$ is
called \textbf{minimal} if $n = \ell(\delta(x_0, x_n))$.

A building is called \textbf{thin} (resp. \textbf{thick}) if all
of its panels have cardinality~2 (resp. at least~3). Any thin
building of type $(W,S)$ is isomorphic to the canonical building
$(W,S)$.

\subsubsection{Projections and convexity in arbitrary buildings}\label{sect:bldg:proj}

A fundamental property of buildings, besides the existence of
apartments, is the exitence of projections onto residues. We
review here the main properties of projections in arbitratry
buildings. The notion of a projection can be slightly refined in
the case of twin buildings; we will come back to this refinement
in Section \ref{sect:twbldg:proj} below.

Let $(W,S)$ be a Coxeter system. We recall from \cite{Bo81} that
if $J, K \subseteq S$ and $w \in W$, then there is a unique
element of minimal length in the double coset $W_J w W_K$.

Let $(\Delta, \delta)$ be a building of type $(W, S)$ and let
$R_J, R_K$ be residues of $\Delta$ of respective type $J$ and $K$.
Then the set of all $\delta(c, d)$ for $c \in R_J$ and $d \in R_K$
is a double coset $W_J w W_K$. Its minimal element is denoted by
$\delta(R_J, R_K)$.

 The set
$$\proj_{R_J}(R_K):= \{c \in R_J | \exists d \in R_K
\text{ such that } \delta(c, d) = \delta(R_J, R_K)\}$$ is called
the \textbf{projection of $R_K$ on $R_J$}. It is a residue of type
$J \cap wK w^{-1}$, where $w:=\delta(R_J, R_K)$; in particular, it
is a spherical residue whenever $J$ or $K$ is spherical. Moreover,
we have
$$\proj_{R_J}(R_K) = \{\proj_{R_J}(c) | c \in
R_K\},$$ where we have written $\proj_{R_J}(c)$ for
$\proj_{R_J}(\{c\})$.

If $c$ is a chamber and $R$ a residue, then $\proj_R(c)$ is a gate
of $c$ to $R$. This means that for any $x \in R$ there exists a
minimal gallery joining $c$ to $x$ via $\proj_R(c)$. The chamber
$\proj_R(c)$ is the unique chamber of $R$ at minimal numerical
distance from $c$.

A set $\mathcal{X}$ of chambers in a building $\Delta$ is called
\textbf{convex} if the following property holds: \emph{given
chambers $x, x' \in \mathcal{X}$ and a spherical residue $R$
containing $x$, then $\proj_R(x') \in \mathcal{X}$.} For example,
apartments and residues are convex sets of chambers.

\subsubsection{Spherical residues and spherical
buildings}\label{sect:bldg:sph}

A building $(\Delta, \delta)$ of type $(W,S)$ is called
\textbf{spherical} if $W$ is finite. In that case, there exists a
unique element $w_0$ of maximal length in $W$. Two chambers $x, y
\in \Delta$ are called \textbf{opposite} if $\delta(x, y) = w_0$.
Two residues $R_J$ and $R_K$ of $\Delta$ of type $J$ and $K$
respectively are called \textbf{opposite} if they contain opposite
chambers and if $J = w_0 K w_0^{-1}$.

A residue $R$ of type $J$ in a building of arbitrary type $(W,S)$
is called \textbf{spherical} if $J$ is a spherical subset of $S$.
Thus $R$ is a spherical building and it makes sense to talk about
opposite chambers and opposite residues of $R$.

The following lemma is a useful criterion of sphericity in terms
of projections.
\begin{lem}\label{criterion of spher}
Let $(\Delta, \delta)$ be a building of type $(W,S)$, let $J
\subseteq S$ and let $R$ be a $J$-residue. Then $J$ is spherical
if and only if there exist $x, y \in R$ such that for every $j \in
J$, we have $\proj_{\pi_j}(x) \neq y$ where $\pi_j$ denotes the
$j$-panel containing $y$.
\end{lem}
\begin{proof} This follows from \cite{Ro89}, Theorem (2.16).
\end{proof}

We end this subsection by recalling a celebrated fixed point
theorem for finite groups acting on buildings.

\begin{prop}\label{fixed point}
Any finite group acting on a building of type $(W, S)$, where $S$
is finite, stabilizes a residue of spherical type.
\end{prop}
\begin{proof} See \cite{Da97}, Corollary 11.9.
\end{proof}

\subsection{Twin buildings}

\subsubsection{Definition}\label{sect:twbldg:def}

A \textbf{twinned pair of buildings} or \textbf{twin building} of
type $(W, S)$ is a pair $((\Delta_+, \delta_+),( \Delta_-,
\delta_-))$ of buildings of type $(W,S)$, endowed with a
\textbf{$W$-codistance}
$$\delta^{\ast}~: (\Delta_+ \times \Delta_-) \cup (\Delta_- \times
\Delta_+) \to W$$ satisfying the following axioms, where $\epsilon
\in \{+, -\}$, $x \in \Delta_{\epsilon}$,  $y \in
\Delta_{-\epsilon}$ and $w = \delta^{\ast}(x, y)$:
$$\begin{array}{ll}
\text{\textbf{(Tw1)}} & \cod(y, x) = w^{-1};\\
\text{\textbf{(Tw2)}} & \text{if } z \in \Delta_{-\epsilon} \text{
is such that } \delta_{-\epsilon}(y, z)= s \in S \text{ and }
\ell(ws) < \ell(w), \text{
then } \cod(x, z)=ws;\\
\text{\textbf{(Tw3)}} & \text{if } s \in S, \text{ there exists }
z \in \Delta_{-\epsilon} \text{ such that } \delta_{-\epsilon}(y,
z) = s \text{ and } \cod(x, z)=ws.
\end{array}$$
In the sequel, we will often use the symbol $\Delta_+$ to denote
the building $(\Delta_+, \delta_+)$ as well as its set of chambers
and similarly for $\Delta_-$. The meaning will be clear from the
context.

A residue $R$ of the twin building $\Delta=(\Delta_+,
\Delta_-,\cod)$ is a residue of $\Delta_\eps$, for $\eps = +$ or
$-$ and $\eps$ is called the \textbf{sign} of $R$. Two chambers
$x$ and $y$ of opposite signs are called \textbf{opposite} if
$\cod(x, y)=1$. Two residues are called \textbf{opposite} if they
are of the same type and contain opposite chambers. Given $J
\subseteq S$, then a pair of opposite residues of type $J$ endowed
with the $W$-codistance induced from $\cod$ is itself a twin
building of type $(W_J, J)$.

Notice that we have defined the term \emph{opposite} at two
different places, namely in Section \ref{sect:bldg:sph} above and
here in Section \ref{sect:twbldg:def}. However, this terminology
is standard and coherent. Indeed, the former notion applies to
chambers or residues of the same sign and lying in a common
spherical residue, while the latter applies to chambers or
residues of opposite signs.

An \textbf{automorphism} of $\Delta=(\Delta_+, \Delta_-, \cod)$ is
by definition a pair $\varphi = (\varphi_+, \varphi_-)$ of
permutations of $\Delta_+$ and $\Delta_-$ respectively preserving
the $W$-distances $\delta_+$, $\delta_-$ as well as the
$W$-codistance $\cod$. \textbf{Isomorphisms} of twin buildings are
defined similarly. We recall from \cite{Ti89}, Theorem 1, that if
$\Delta_+$ and $\Delta_-$ are thick, then an automorphism of
$\Delta$ which fixes a pair of opposite chambers $c, c'$ and all
chambers adjacent to $c$ is the identity.

\subsubsection{Reflections and twin apartments}\label{sect:twbldg:apt}

Let $(\Sigma_+, \delta_+)$ and $(\Sigma_-, \delta_-)$ be two
copies of the canonical building $\mathcal{A}(W,S)$ of type
$(W,S)$ (see Section \ref{sect:bldg:apt}). Let $\cod : \Sigma_+
\times \Sigma_- \to W : (x,y) \mapsto x^{-1}y$; this makes sense
since $\Sigma_+=\Sigma_-=W$. Then $\Sigma(W, S) := ((\Sigma_+,
\delta), (\Sigma_-, \delta), \cod)$ is a thin twin building of
type $(W,S)$, namely a twinned pair of thin buildings. It is the
unique thin twin building of that type up to isomorphism.

The group $W$ has a faithful action on $\Sigma(W,S)$ by
automorphisms which is given by left multiplication on $\Sigma_+$
and $\Sigma_-$. Every automorphism of $\Sigma$ is of this form.

A \textbf{reflection} is a non-trivial element of $W$ which
stabilizes a panel of $\Sigma(W,S)$. Conversely, to any panel of
$\Sigma(W,S)$ corresponds a unique reflection of $W$ which
stabilizes it. Moreover, an element of $W$ is a reflection if and
only if it is conjugate to an element of $S$.

Let $\Delta=(\Delta_+, \Delta_-)$ be a twin buildings of type
$(W,S)$. A pair $\Sigma=(\Sigma_+, \Sigma_-)$ of subsets of
$\Delta$ is called a \textbf{twin apartment} if it is isomorphic
to the canonical twin building $\Sigma(W,S)$. Given a twin
apartment $\Sigma=(\Sigma_+, \Sigma_-)$, the restriction of the
opposition relation of $\Delta$ to $\Sigma$ is a one-one
correspondence $\Sigma_+ \leftrightarrow \Sigma_-$. (It
corresponds to the identity $\mathrm{id} : W \to W$ in the
canonical twin building $\Sigma(W,S)$.) We denote it by
$\op_\Sigma$.

It is a fundamental fact that, given any two chambers $x \in
\Delta_\eps$ and $y \in \Delta_{\eps'}$ in a twin building
$\Delta=(\Delta_+, \Delta_-)$ of type $(W,S)$, where $\eps, \eps'
\in \{+,-\}$, there exists a twin apartment $\Sigma=(\Sigma_+,
\Sigma_-)$ such that $x \in \Sigma_\eps$ and $y \in
\Sigma_{\eps'}$ (see \cite{Ab96}, Lemma 2). It is common and handy
to say that $x$ and $y$ are contained in $\Sigma$, and to write
$x, y \in \Sigma$.

\subsubsection{Projections and convexity in twin buildings}\label{sect:twbldg:proj}

Let $(W,S)$ be a Coxeter system and let $J, J \subset S$ be
spherical subsets. Then there is a unique element of maximal
length in $W_J w W_K$ (\cite{Ab96}, Lemma 9).

Let $\Delta= (\Delta_+, \Delta_-, \cod)$ be a twin building of
type $(W,S)$ and let $R_J, R_K$ be residues of $\Delta$ of
respective type $J$ and $K$. Assume moreover that $J$ and $K$ are
spherical and that $R_J$ and $R_K$ have opposite signs. Then the
set of all $\cod(c, d)$ for $c \in R_J$ and $d \in R_K$ is a
double coset $W_J w W_K$. Its maximal element is denoted by
$\cod(R_J, R_K)$. The set
$$\proj_{R_J}(R_K):= \{c \in R_J | \exists d \in R_K
\text{ such that } \cod(c, d) = \cod(R_J, R_K)\}$$ is called the
\textbf{projection of $R_K$ on $R_J$}. It is a residue of type
$w^0_J(J \cap w K w^{-1})w^0_J$, where $w:=\cod(R_J, R_K)$ and
$w^0_J$ denotes the maximal element of $W_K$ (\cite{Ab96}, Lemma
10). Moreover, we have
$$\proj_{R_J}(R_K) = \{\proj_{R_J}(c) | c \in
R_K\}.$$

A set $\mathcal{X}$ of chambers in a twin building is called
\textbf{convex} if the following condition holds: \emph{given $x,
x' \in \mathcal{X}$ and a spherical residue $R$ containing $x'$,
then $\proj_R(x) \in \mathcal{X}$.} For example, twin apartments
are convex sets of chambers in twin buildings. Actually, the
convex hull of any pair of opposite chambers is a twin apartment
containing them. This implies that two opposite chambers lie in a
unique common twin apartment.

Notice that we have defined the terms \emph{projections} and
\emph{convexity} at two different places, namely in Section
\ref{sect:bldg:proj} above and here in Section
\ref{sect:twbldg:proj}. The point is that the notion of
projections in twin buildings is a generalization of the standard
notion of projections in arbitrary buildings. There will be no
confusion between both. Indeed, the meaning of the symbol
$\proj_R(x)$ in the context of twin buildings depends on the
respective signs of the residue $R$ and the chamber $x$.

We end this subsection with a result is often useful to compute
projections between residues of opposite signs using twin
apartments.

\begin{lem}\label{Abr prop4}
Let $\Delta=(\Delta_+, \Delta_-, \cod)$ be a twin building and for
each sign $\eps$ let $R_\eps$ be a spherical residue of
$\Delta_\eps$. Let $\Sigma$ is a twin apartment containing $R_+$
and $R_-$ and let $\eps \in \{+,-\}$. Let $R'_\eps$ be the residue
of $\Sigma$ opposite $R_{-\eps}$. Then the residues
$\proj_{R_\eps}(R_{-\eps})$ and $\proj_{R_\eps}(R'_{\eps})$ are
contained in $\Sigma$ and opposite in $R_\eps$ (see Section
\ref{sect:bldg:sph}).
\end{lem}
\begin{proof} This is Proposition 4 in \cite{Ab96}.
\end{proof}

In brief, the statement of this lemma may be written as
$$\proj_{R_\eps}(R_{-\eps}) = \op_{\Sigma \cap
R_\eps} (\proj_{R_\eps}(\op_\Sigma(R_{-\eps}))).$$

\subsubsection{Parallelism}\label{sect:twbldg:parallel}

Let $\Delta = (\Delta_+, \Delta_-, \cod)$ be a twin building of
type $(W,S)$. Two residues $R_J, R_K$ of $\Delta$ (assumed to be
spherical if they have opposite signs) are called
\textbf{parallel} if $\proj_{R_J}(R_K)=R_J$ and $\proj_{R_K}(R_J)
=R_K$.

It follows from the definitions that $\proj_{R_J}(R_K)$ and
$\proj_{R_K}(R_J)$ are always parallel.

Although parallel residues need not have the same type, they are
nevertheless always `almost isometric' in the following sense.
\begin{lem}\label{isom parallel}
Let $R_J$ (resp. $R_K$) be a residue of spherical type $J$ (resp.
$K$) and sign $\eps_J$ (resp. $\eps_K$). Assume that $R_J$ and
$R_K$ are parallel. Then there exists an isomorphism $\eta : W_J
\to W_K$ with $\eta(J) = K$ such that
$$\delta_{\eps_J}(\proj_{R_J}(x),
\proj_{R_J}(y))=\eta (\delta_{\eps_K}(x, y))$$ for all $x, y \in
R_K$. In particular, if $x$ and $y$ are opposite in $R_K$, then so
are $\proj_{R_J}(x)$ and $\proj_{R_K}(y)$ in $R_J$.
\end{lem}

The proof of Lemma \ref{isom parallel} is in the same spirit as
the proof of Proposition 5.15 in \cite{We03} and is omitted here.
We only mention that the isomorphism $\eta$ of the lemma is
actually induced by the conjugation by $\delta_{\eps_J}(R_J, R_K)$
if $\eps_J = \eps_K$ and by $w_J \cod(R_J, R_K)$ if $\eps_J =
-\eps_K$. However, we do not need this fact here.

\begin{lem}\label{iff opposite}
The spherical residues $R_J$ and $R_K$ of $\Delta$ are opposite if
and only if $\proj_{R_J}(R_K)$ and $\proj_{R_K}(R_J)$ are
opposite.
\end{lem}
\begin{proof} Since opposite spherical residues are parallel, the
implication `$\Rightarrow$' is obvious. The other implication
follows from an easy computation using Lemma 9 of \cite{Ab96}.
\end{proof}

Next, we give a rule for the composition of projections.
\begin{lem}\label{compose proj}
As before, $\Delta$ is a (possibly twin) building of type $(W,
S)$. Let $R_I, R_J, R_K$ be residues of type $I, J, K$
respectively and assume that $R_I \subseteq R_J$. Moreover, if
$\Delta$ is a twin building and if $R_J$ and $R_K$ have opposite
signs, then we also assume that $I, J, K$ are spherical. Then we
have
$$\proj_{R_I}(R_K) =
\proj_{R_I}(\proj_{R_J}(R_K)).$$
\end{lem}
\begin{proof} It suffices to prove the statement when the residue $R_K$ is
reduced to a single chamber, say $c$ (or, in other words, when $K
= \emptyset$). If $R_J$ and $c$ have the same sign, the result
follows from the fact that $\proj_{R_J}(c)$ is a gate of $c$ to
$R_J$. If they have opposite signs, we may reduce ourselves to the
preceding case in view of Lemma \ref{Abr prop4}.
\end{proof}

The following lemma characterizes the parallelism of spherical
residues in thin buildings.

\begin{prop}\label{parallelism}
Let $(\Sigma, \delta)$ be the thin building of type $(W,S)$. Let
$J, K$ be spherical subsets of $S$ and let $R_J, R_K$ be residues
of type $J, K$ respectively. Then the following statements are
equivalent~:\\
\begin{tabularx}{\linewidth}{rX}
(i) & $R_J$ and $R_K$ are parallel;\\
(ii) & a reflection of $\Sigma$ stabilizes $R_J$ if and only if
it stabilizes $R_K$;\\
(iii) & there exist two sequences $R_J = R_0, R_1, \dots, R_n =
R_K$ and $T_1, \dots, T_n$ of residues of spherical type such that
for each $1 \leq i \leq n$ the rank of $T_i$ is equal to $1 +
\mathrm{rank}(R_J)$, the residues $R_{i-1}$, $R_i$ are contained
and opposite in $T_i$ and moreover, we have
$\proj_{T_i}(R_J)=R_{i-1}$ and $\proj_{T_i}(R_K)=R_i$.
\end{tabularx}
\end{prop}
\begin{proof} The equivalence (i) $\Leftrightarrow$ (ii) is easy. The
implication (iii) $\Rightarrow$ (i) follows from an obvious
induction on $n$ using the fact that opposite spherical residues
are parallel. It remains to prove (i) $\Rightarrow$ (iii). Let $s
\in S$ such that $\ell(s\delta(R_J, R_K)) < \ell(\delta(R_J,
R_K))$. Clearly $s \not \in J$. Let $x \in R_J$ and set
$T_1:=\mathrm{Res}_{J \cup \{s\}}(x)$ and $R_1 :=
\proj_{T_1}(R_K)$. By definition of $T_1$ we have $R_J \cap R_1 =
\emptyset$ and so $R_1$ is properly contained in $T_1$. By Lemma
\ref{compose proj} we have $\proj_{R_J}(R_1)= R_J$. Therefore
$R_1$ and $R_J$ have the same rank and so they are parallel.

Let $x' \in R_1$ such that $\proj_{R_J}(x')=x$ and choose $y$
opposite to $x'$ in $R_1$ (this makes sense since $R_1$, being the
image of $R_K$ under a projection, is spherical). Let now $\pi$ be
a panel containing $x$ and contained in $T_1$. If the type of
$\pi$ is an element of $J$ then $\proj_\pi(y) \neq x$ by Lemma
\ref{criterion of spher} and Lemma \ref{isom parallel}. If the
type of $\pi$ is $s$ then the same inequality is still true by the
definition of $s$ and using Lemma \ref{compose proj}. Therefore,
$T_1$ is spherical by Lemma \ref{criterion of spher}. Since
$\delta(R_1, R_K)$ is shorter than $\delta(R_J, R_K)$ by
construction, the desired result follows from an immediate
induction.
\end{proof}

\begin{cor}\label{proj max sph}
Let $\Delta$ be a (possibly twin) building of type $(W,S)$ and let
$R_K$ be a spherical residue which is maximal with respect to that
property. Then, for any residue $R_J$ (assumed to be spherical if
$\Delta$ is a twin building and if $R_J$  and $R_K$ have opposite
signs) the projection of $R_J$ on $R_K$ is properly contained in
$R_K$ unless $R_J$ and $R_K$ are equal or opposite.
\end{cor}
\begin{proof} Since $\proj_{R_K}(R_J)$ and
$\proj_{R_J}(R_K)$ are parallel, the result clearly follows from
the previous proposition, using also Lemma \ref{Abr prop4} if
$R_J$ and $R_K$ have opposite signs.
\end{proof}

\begin{cor}\label{equiv relation}
Let $\Delta$ be a twin building and let $\Sigma$ be a twin
apartment of $\Delta$. Then the parallelism is an equivalence
relation on the set of spherical residues of $\Sigma$.
\end{cor}
\begin{proof} This follows from Lemma \ref{Abr prop4} and Proposition
\ref{parallelism}.
\end{proof}

\subsubsection{Twin roots}\label{sect:twbldg:twroot}

Let $\Delta= (\Delta_+, \Delta_-, \cod)$ be a twin building of
type $(W, S)$.

A \textbf{twin root} of $\Delta$ is the convex hull of a pair of
chambers `at codistance $1$', namely a pair $\{x, y\}$ such that
$s := \cod(x, y) \in S$. Let $\pi$ be the $s$-panel containing
$x$. Then any chamber $x' \in \pi \backslash \{x\}$ is opposite
$y$ and determines therefore a twin apartment $\Sigma$ which
contains $\phi$ and $x'$ (see Section \ref{sect:twbldg:proj}). We
say that $\phi$ is a twin root of $\Sigma$. The complement of
$\phi$ in $\Sigma$ is also a twin root; it is actually the convex
hull of $x'$ and $\proj_{\pi}(x')$. This twin root is said to be
\textbf{opposite} to $\phi$ in $\Sigma$ and is denoted by $-\phi$
although its definition depends on $\Sigma$. A residue $R$ of
$\Delta$ is said to be \textbf{in the interior} of $\phi$ if it is
contained in $\Sigma$ and if $R \cap \Sigma$ is contained in
$\phi$. If $R \cap \phi$ and $R \cap (-\phi)$ are both nonempty,
then $R$ is said to be \textbf{on the boundary} of $\phi$.

\subsection{From groups to buildings: twin root data}\label{sect:TRD}

\subsubsection{Definition}\label{sect:TRD:defPrinc}
Let $\Sigma = \Sigma(W, S)$ be the canonical twin building of type
$(W, S)$ (see Section \ref{sect:twbldg:apt}) and let $\Phi(W,S)$
be the set of all its twin roots. We have already mentioned the
action of $W$ on $\Sigma$ (see Section \ref{sect:twbldg:apt}).
Given a twin root $\phi \in \Phi(W ,S)$, then all panels on the
boundary of $\phi$ correspond to the same reflection of $W$. This
reflection is denoted by $s_\phi$ and it permutes $\phi$ and
$-\phi$. A pair $\{\phi, \psi\}$ of twin roots of $\Sigma=
(\Sigma_+, \Sigma_-)$ is said to be \textbf{prenilpotent} if $\phi
\cap \psi \cap \Sigma_+$ and $(-\phi) \cap (-\psi) \cap \Sigma_+$
are both nonempty; in that case, we denote by $[\phi, \psi]$ the
set of all twin roots $\alpha$ of $\Sigma$ such that $\alpha
\supseteq \phi \cap \psi$ and $-\alpha \supseteq (-\phi) \cap
(-\psi)$.

A \textbf{twin root datum} of type $(W,S)$ is a system
$\mathcal{D}:=(G, (U_\phi)_{\phi \in \Phi(W,S)})$ consisting of a
group $G$ and a family of subgroups $U_{\phi}$ which satisfies the
following axioms, where $H$ and $U(c)$ denote respectively the
intersection of the normalizers of all $U_{\phi}$'s and the
subgroup of $G$ generated by the
$U_{\phi}$'s such that $\phi$ contains the chamber $c$ of $\Sigma$~:\vspace{.3cm}\\
\begin{tabularx}{\linewidth}{lX}
\textbf{(TRD0)} & $U_{\phi} \not = 1 $ for all $\phi \in \Phi(W,S)$;\\
\textbf{(TRD1)} & if $\{\phi, \psi\}$ is a prenilpotent pair of
distinct twin roots, the commutator $[U_{\phi}, U_{\psi}]$ is
contained in the group generated by all $U_{\gamma}$'s with
$\gamma \in
[\phi, \psi] \backslash \{\phi, \psi\}$;\\
\textbf{(TRD2)} & if $\phi \in \Phi(W, S)$ and $u \in
U_{\phi}\backslash \{1\}$, there exists elements $u', u''$ of
$U_{-\phi}$ such that the product $\mu(u)=u'uu''$
conjugates $U_{\psi}$ onto $U_{s_{\phi}(\psi)}$ for each $\psi \in \Phi(W,S)$; \\
\textbf{(TRD3)} & if $\phi \in \Phi(W, S)$ and $c$ is a chamber of
$\Sigma$ which is not contained in $\phi$, then $U_{\phi}$ is not
contained in $U(c)$;\\
\textbf{(TRD4)} & the group $G$ is generated by $H$ and the
$U_{\phi}$'s.
\end{tabularx}
The group $G$ is sometimes denoted by $G^\mathcal{D}$.

\subsubsection{Isomorphisms of twin root data}\label{sect:TRD:def}

Let $\mathcal{D}:=(G, (U_\phi)_{\phi \in \Phi(W,S)})$ and
$\mathcal{D}':=(G', (U'_\phi)_{\phi \in \Phi(W',S')})$ be twin
root data. Let $S = S_1 \cup \dots \cup S_n$ be the finest
partition of $S$ such that $[S_i, S_j]=1$ whenever $1 \leq i < j
\leq n$. Then $\mathcal{D}$ and $\mathcal{D}'$ are called
\textbf{isomorphic} if there exist an isomorphism $\varphi : G \to
G'$, an isomorphism $\pi : W \to W'$ with $\pi(S) = S'$, an
element $x \in G'$ and a sign $\eps_i$ for each $1 \leq i \leq n$
such that
\begin{equation}\label{eq1}
x \varphi( U_{\phi} ) x^{-1} = U'_{\eps_i \pi(\phi)} \qquad \qquad
\text{for every twin root $\phi$ with $s_\phi \in W_{S_i}$,}
\end{equation}
where we denote by $\pi$ the obvious bijection $\Phi(W, S) \to
\Phi(W', S)$ induced by $\pi : W \to W'$. Thus, if $(W,S)$ is
irreducible, then either $ x \varphi( U_{\phi} ) x^{-1} =
U'_{\pi(\phi)}$ or $x \varphi( U_{\phi} ) x^{-1} =
U'_{-\pi(\phi)}$ for all $\phi \in \Phi(W,S)$.

When (\ref{eq1}) holds, we say that the isomorphism $\varphi$
\textbf{induces an isomorphism} of $\mathcal{D}$ to
$\mathcal{D}'$. In particular, this means, the isomorphism
$\varphi$ maps the union of conjugacy classes
$$\{gU_+g^{-1} | g \in G\} \cup \{gU_-g^{-1} | g \in G\}$$
to
$$\{gU'_+g^{-1} | g \in G'\} \cup \{gU'_-g^{-1} | g \in G\},$$
where $U_+$ (resp. $U_-$, $U'_+$, $U'_-$) denotes $U(c)$ (resp.
$U(\op_\Sigma(c))$, $U(c')$, $U(\op_{\Sigma'}(c'))$) for some $c
\in \Sigma=\Sigma(W,S)$ and some $c' \in \Sigma'=\Sigma(W',S')$.

A crucial fact on isomorphisms between twin root data we will need
later is the following.

\begin{prop}\label{thm:KM1}
Let $\mathcal{D}:=(G, (U_\phi)_{\phi \in \Phi(W,S)})$ and
$\mathcal{D}':=(G', (U'_\phi)_{\phi \in \Phi(W',S')})$ be twin
root data with $S$ and $S'$ finite and let $\varphi: G \to G'$ be
an isomorphism. Assume there exists $g \in G'$ such that
$$\{\varphi(U_\phi) | \; \phi \in \Phi(W,S)\} = \{g U'_\phi g^{-1}
| \; \phi \in \Phi(W',S')\}.$$ Then $\varphi$ induces an
isomorphism of $\mathcal{D}$ to $\mathcal{D'}$.
\end{prop}
\begin{proof}
This is Theorem (2.5) in \cite{CM03b}.
\end{proof}

\subsubsection{Twin buildings from twin root
data}\label{sect:TRD:bldg}

Let  $\mathcal{D}:=(G, (U_\phi)_{\phi \in \Phi(W,S)})$ be a twin
root datum of type $(W, S)$. Let $H$ be the intersection of the
normalizers of all $U_{\phi}$'s and let $N$ be the subgroup of $G$
generated by $H$ together with all $\mu(u)$ such that $u \in
U_\phi \backslash \{1\}$, where $\mu(u)$ is as in (TRD2). Let $c$
be a chamber of $\Sigma = \Sigma(W, S)$ of positive sign, let $c'
:= \op_\Sigma(c)$ and let $B_+ := H.U(c)$ and $B_-:=H. U(c')$.

We recall from \cite{Ti92}, Proposition 4, that $(G, B_+, N)$ and
$(G, B_-, N)$ are both $BN$-pairs of type $(W,S)$. Thus, we have
corresponding Bruhat decompositions of $G$:
$$G = \coprod_{w \in W} B_+ w B_+ \qquad \qquad \text{and} \qquad \qquad
G = \coprod_{w \in W} B_- w B_-.$$ For each $\eps \in \{+,-\}$,
the set $\Delta_\eps := G/B_\eps$ endowed with the map
$\delta_\eps: \Delta_\eps \times \Delta_\eps \to W$ by
$$\delta_\eps(gB_\eps, hB_\eps)=w \Leftrightarrow B_\eps g^{-1} h
B_\eps = B_\eps w B_\eps,$$
has a canonical structure of a thick
building of type $(W,S)$.

The twin root datum axioms imply that $G$ also admits Birkhoff
decompositions (see Lemma 1 in \cite{Ab96}):
$$G = \coprod_{w \in W} B_\eps w B_{-\eps}$$
for each $\eps \in \{+,-\}$. The pair $((\Delta_+,\delta_+),
(\Delta_-, \delta_-))$ of buildings admits a natural twinning by
means of the $W$-codistance $\cod$ defined by
$$\cod(gB_\eps, h
B_{-\eps})=w \Leftrightarrow B_\eps g^{-1} h B_{-\eps} = B_\eps w
B_{-\eps}$$ for each $\eps \in \{+,-\}$. The triple $\Delta :=
((\Delta_+, \delta_+), (\Delta_-, \delta_-), \cod)$ is a twin
building of type $(W, S)$.

We may and shall identify the chamber $c$ (resp. $c'$) of
$\Sigma(W, S)$ with the chamber $B_+$ of $\Delta_+ = G/B_+$ (resp.
$B_-$ of $\Delta_-$). We also identify $\Sigma(W, S)$ with the
unique twin apartment of $\Delta$ containing $c$ and $c'$; this
twin apartment is denoted by $\Sigma$ and is called the
\textbf{fundamental twin apartment} of $\Delta$ (with respect to
the twin root datum $\mathcal{D}$).

The diagonal action of $G$ on $\Delta_+ \times \Delta_-$ by left
multiplication is transitive on pairs of opposite chambers and,
hence, on twin apartments.

\subsubsection{Parabolic subgroups and root subgroups}\label{sect:TRD:parab}

We keep the notation of the previous subsection. We recall from
the theory of $BN$-pairs (see \cite{Bo81}, Chapter IV) that a
subgroup $P$ of $G$ containing $B_\eps$ or any of its conjugates
is called a \textbf{parabolic subgroup} of sign $\eps$, where
$\eps \in \{+,-\}$. If $P$ contains $B_\eps$, then there exists $J
\subseteq S$ such that $P$ has a Bruhat decomposition
$$P= \coprod_{w \in W_J} B_\eps w B_\eps;$$ the set $J$ is called
the \textbf{type} of the parabolic subgroup $P$. If $J$ is
spherical, then $P$ is said to be \textbf{of finite type} (or of
spherical type). A minimal parabolic subgroup (i.e. a parabolic
subgroup of type $\emptyset$) such as $B_+$ or $B_-$ is called a
\textbf{Borel subgroup}.

For $\eps \in \{+,-\}$, let $P^J_\eps$ be the parabolic subgroup
of type $J$ containing $B_\eps$. The geometric meaning of the
groups $B_+$, $B_-$, $P^J_+$, $P^J_-$, $H$ and $N$ is as follows:
$$B_+=\Stab_G(c), \quad B_-=\Stab_G(c'),
\quad P^J_+= \Stab_G(\Res_J(c)), \quad P^J_-=
\Stab_G(\Res_J(c'))$$ and
$$N=\Stab_G(\Sigma),
\qquad H = B_+ \cap N = B_- \cap N = \Fix_G(\Sigma).$$

Given a twin root $\phi$ of $\Sigma$, then $U_\phi$ fixes
chamberwise any panel in the interior of $\phi$ and is sharply
transitive on the set of twin apartments containing $\phi$.
Moreover, it follows then from the axioms that for each $g \in G$
and each twin root $\phi$ of $\Sigma$, the group $U_{g(\phi)}:=g
U_\phi g^{-1}$ depends only on the twin root $g(\phi)$ and not on
the choice of $g$ and $\phi$. Hence, for every twin root $\psi$ of
$\Delta$, there is a well defined group $U_\psi$ which fixes
chamberwise any panel in the interior of $\psi$. The group
$U_\psi$ is sharply transitive on the set of twin apartments
containing $\psi$; it is called the \textbf{root subgroup}
associated with the twin root $\psi$.

\subsubsection{The Conditions (P1)--(P3) and a technical
lemma}\label{sect:TRD:condP}

Let $\mathcal{D}:=(G, (U_{\phi})_{\phi \in \Phi(W,S)})$ be a twin
root datum. For each $\phi \in \Phi(W, S)$, we set $L_\phi :=
\langle U_\phi \cup U_{-\phi} \rangle$ and $H_\phi :=
N_{L_\phi}(U_\phi) \cap N_{L_\phi}(U_{-\phi})$. The group $H_\phi$
acts on the conjugacy class $\mathcal{C}_\phi$ of $U_\phi$ in
$L_\phi$. We shall be interested in the following three conditions
(see Theorem \ref{main} below)~:\\
\begin{tabularx}{\linewidth}{rX}
\textbf{(P1)} & for every $\phi \in \Phi(W, S)$, the group
$U_\phi$ is nilpoptent;\\
\textbf{(P2)} & for every $\phi \in \Phi(W, S)$, the group
$L_\phi$ is perfect;\\
\textbf{(P3)} & for every $\phi \in \Phi(W, S)$, the groups
$U_\phi$ and $U_{-\phi}$ are the only fixed points of $H_\phi$ in
$\mathcal{C}_\phi$.
\end{tabularx}

The following lemma gives the geometric interpretation of
Condition (P3).

\begin{lem}\label{H not trivial}
Let $\mathcal{D}= (G, (U_\phi)_{\phi \in \Phi(W, S)})$ be a twin
root datum of type $(W, S)$ which satisfies Condition (P3). Let
$\Delta$ be the twin building associated with $\mathcal{D}$ and
let $H$ be as above (see Section \ref{sect:TRD:bldg}). Then $H$
fixes no chamber outside the fundamental twin apartment of
$\Delta$.
\end{lem}
\begin{proof} Let $\Sigma$ be the fundamental twin apartment of $\Delta$.
Let $\phi \in \Phi$ and let $H_\phi$ be the intersection of the
normalizers of $U_\phi$ and $U_{-\phi}$ in $\langle U_\phi \cup
U_{-\phi}\rangle$. The group $H_\phi$ fixes $\Sigma$ chamberwise
and is therefore contained in $H$. Let $\pi$ be any panel on the
boundary of $\phi$. Then $\pi$ is stabilized by $U_\phi$ and
$U_{-\phi}$. The condition (P3) means precisely that the only
fixed chambers of $H_\phi$ in $\pi$ are the two elements of $\pi
\cap \Sigma$. This implies that for any panel $\pi$ of $\Sigma$,
the group $H$ fixes no chamber in $\pi \backslash \Sigma$. Now, an
easy induction on the gallery distance from an arbitrary chamber
of $\Delta$ to $\Sigma$ finishes the proof.
\end{proof}

\subsubsection{Twin root data over fields}\label{sect:TRD:fields}

Let $\mathcal{D}= (G, (U_\phi)_{\phi \in \Phi(W, S)})$ be a twin
root datum of type $(W, S)$, let $H \leq G$ be as in the previous
subsection and for each $\phi \in \Phi(W, S)$, let $\K_\phi$ be a
field. We recall from \cite{CM03b} that the twin root datum
$\mathcal{D}$ is called \textbf{locally split over the fields
$(\K_\phi)_{\phi \in \Phi(W, S)}$} if $H$ is abelian and if for
each $\phi \in \Phi(W, S)$, the twin root datum $\mathcal{D}_\phi
:= (H\langle U_\phi \cup U_{-\phi}\rangle, \{U_\phi, U_{-\phi}\})$
is isomorphic to the natural twin root datum of $SL_2(\K_\phi)$ or
$PSL_2(\K_\phi)$. Of course, the natural twin root datum
associated to a (split) Kac-Moody group over a field $\K$ is
locally split over $\K$.

Notice that if $\mathcal{D}$ is locally split over the fields
$(\K_\phi)_{\phi \in \Phi(W, S)}$, then $\mathcal{D}$ satisfies
Condition~(P1). Moreover, if $\K_\phi$ has cardinality at least 4
for every $\phi \in \Phi(W, S)$, then Conditions~(P2) and (P3) are
also satisfied.

\section{Levi decomposition in twin root data}\label{Levi sect}

The purpose of this section is to obtain a Levi decomposition for
intersections of finite type parabolic subgroups of opposite signs
in a group with twin root datum (see Proposition~\ref{Levi for
parallel}). In the language of buildings, this group is the
stabilizer of a pair of spherical residues of opposite signs.

\subsection{Levi decomposition of parabolic subgroups}

\subsubsection{The setting}
\label{sect:Levi:notation}

 Let $\mathcal{D}:=(G, (U_\alpha)_{\alpha \in
\Phi(W,S)})$ be a twin root datum, let $\Delta=((\Delta_+,
\delta_+), (\Delta_-, \delta_-), \cod)$ be the corresponding twin
building and let $\Sigma_0$ be the fundamental twin apartment (see
Section \ref{sect:TRD:bldg}).

Let $\Sigma$ be any twin apartment of $\Delta$. Let $c \in \Sigma$
be a chamber and let $R$ be a spherical residue of $\Sigma$ (i.e.
$R \cap \Sigma \neq \emptyset$). Let $\Phi^\Sigma$ the set of all
twin roots of $\Sigma$ and let $\Phi^\Sigma(R)$ be the set of twin
roots $\beta$ of $\Sigma$ such that $R \cap \beta$ and $R \cap
(-\beta)$ are both nonempty, which means precisely that the
reflection $s_\beta$ stabilizes $R$ (see
Sections~\ref{sect:twbldg:apt} and~\ref{sect:twbldg:twroot}). We
set
$$U^\Sigma(c):= \la U_\phi | \; \phi \in \Phi^\Sigma \ra,
\qquad U^\Sigma(R) := \bigcap_{x \in \Sigma \cap R} U^\Sigma(x)$$
and
$$L^\Sigma(R) := \Fix_G(\Sigma). \la U_\phi | \; \phi \in
\Phi^\Sigma(R) \ra.$$ (See Section \ref{sect:TRD:parab} for the
definition of $U_\phi$.)

We will see in Proposition~\ref{Levi} that $U^\Sigma(c)$ and
$U^\Sigma(R)$ are actually independent of $\Sigma$. They will be
denoted by $U(c)$ and $U(R)$ respectively.

Notice that $L^\Sigma(R) = L^\Sigma(\op_\Sigma(R))$ since
$\Phi^\Sigma(R) = \Phi^\Sigma(\op_\Sigma(R))$. Moreover, if the
residue $R$ is reduced to a chamber $c$ (i.e. if $R$ is of type
$\emptyset$) then $L^\Sigma(c)=\Fix_G(\Sigma)$, which is
$G$-conjugate to the subgroup $H$ of Section~\ref{sect:TRD:bldg}
(see Section \ref{sect:TRD:parab}).

\subsubsection{Standard Levi decomposition:
Levi factor and unipotent radical}\label{sect:Levi:standard}

The following result is the standard Levi decomposition of finite
type parabolic subgroups of a group with a twin root datum. We
state it in the language of buildings.

\begin{prop}\label{Levi}
We have $\mathrm{Stab}_G(R) = L^\Sigma(R) \ltimes U^\Sigma(R)$.
Moreover, $U(R)$ is sharply transitive on the set of residues
which are opposite $R$ in $\Delta(G)$ and
$$L^\Sigma(R) =
\Stab_G(R) \cap \Stab_G(\op_\Sigma(R)).$$ In particular, the
subgroup $U^\Sigma(R)$ is independent of $\Sigma$ and will be
denoted by $U(R)$.
\end{prop}
\begin{proof} This follows from the theorem of Section 6.2.2 in
\cite{Re99}.
\end{proof}

The group $U(R)$ is called the \textbf{unipotent radical} of the
parabolic subgroup $\mathrm{Stab}_G(R)$ with respect to the twin
root datum $\mathcal{D}$ and the group $L^\Sigma(R)$ is called a
\textbf{Levi factor}.

\subsection{Levi decomposition of parabolic intersections}

\subsubsection{More definitions and notation}\label{sect:doubleLevi:notation}
We keep the notation of Section~\ref{sect:Levi:notation}.

Let $\Sigma$ be a twin apartment of $\Delta$ and for each $\eps
\in \{+,-\}$, let $R_\eps$ be a residue of $\Sigma$ of sign
$\eps$. We set
$$U^\Sigma(R_+, R_-):= \la U_\beta | \;
\beta \in \Phi^\Sigma \text{ and } R_\eps \cap \Sigma \subset
\beta \ra,$$ and, for $\eps \in \{+, -\}$,
$$
\widetilde{U}^\Sigma(R_\eps, R_{-\eps}):= \la U_\beta | \; \beta
\in \Phi^\Sigma(R_\eps) \text{ and } \proj_{R_\eps}(R_{-\eps})
\subset \beta \ra$$ (see Sections~\ref{sect:twbldg:twroot}
and~\ref{sect:TRD:parab}). Notice that $U^\Sigma(R_+, R_-) =
U^\Sigma(R_-, R_+)$ while $\widetilde{U}^\Sigma(R_+, R_-)$ need
not equal $\widetilde{U}^\Sigma(R_-, R_+)$. If $R_+$ (resp. $R_-$)
is reduced to a chamber $x_+$ (resp. $x_-$), we have
$\widetilde{U}^\Sigma(R_+, R_-) = \widetilde{U}^\Sigma(R_-, R_+) =
\{1\}$. This remains true in the case where $R_+$ and $R_-$ are
parallel (see Section~\ref{sect:twbldg:parallel}).

See Remark~\ref{rem:Utilde} below for an interpretation of the
group $\widetilde{U}^\Sigma(R_\eps, R_{-\eps})$.

\subsubsection{The intersection of a parabolic subgroup and
unipotent radical}

In order to obtain the Levi decomposition of the group
$\Stab_G(\{R_+, R_-\})= \Stab_G(R_+) \cap \Stab_G(R_-)$ , we need
a decomposition result for $\Stab_G(R_+) \cap U(R_-)$. We start
with the case where both residues are reduced to single chambers.

\begin{lem}\label{U_w} Let $x_+ \in \Delta_+$ and $x_- \in
\Delta_-$ be chambers of $\Delta$. Let $\Sigma$ be a twin
apartment containing them both. Let $\eps \in \{+,-\}$, let
$y_\eps := \op_\Sigma(x_{-\eps})$ and let $x_\eps =x_0, x_1,
\dots, x_n=y_\eps$ be a minimal gallery joining $x_\eps$ to
$y_\eps$. For each $1 \leq i \leq n$, let $\beta_i$ be the twin
root of $\Sigma$ containing $x_{i-1}$ but not $x_i$. Then we have
$$U(x_\eps) \cap \mathrm{Stab}_G(x_{-\eps}) = U^\Sigma(x_+, x_-) =
U_{\beta_1}.U_{\beta_2} \dots U_{\beta_n}.$$ In particular, the
product $U_{\beta_1}.U_{\beta_2} \dots U_{\beta_n}$ is a group
which coincides with $U^\Sigma(x_+, x_-)$, and the latter does not
depend on the twin apartment $\Sigma$. We will denote it by
$U(x_+, x_-)$.
\end{lem}
\begin{proof} This follows from Lemma 1.5.2(iii) and Theorem 3.5.4 in
\cite{Re99}.
\end{proof}

The following lemma generalizes Lemma \ref{U_w} to the case of
spherical residues of higher rank.

\begin{lem}\label{U(R_+, R_-)}
Let $R_+ \subseteq \Delta_+$ and $R_- \subseteq \Delta_-$ be
spherical residues of $\Delta$. Let $\Sigma$ be a twin apartment
intersecting them both. Then, for each $\eps \in \{+, -\}$, we
have
$$U(R_{-\eps}) \cap \mathrm{Stab}_G(R_{\eps})=\widetilde{U}^\Sigma(R_\eps, R_{-\eps}) .U^\Sigma(R_+, R_-).$$

In particular, if $R_+$ and $R_-$ are parallel, then
$$U(R_+) \cap \mathrm{Stab}_G(R_-) = \mathrm{Stab}_G(R_+) \cap U(R_-) =U^\Sigma(R_+, R_-).$$
\end{lem}
\begin{proof}
The inclusion `$\supseteq$' of the first part is clear.

Let $R'_\eps := \op_\Sigma(R_{-\eps})$. Let us choose $z, z' \in
\proj_{R'_\eps}(R_\eps) \cap \Sigma$ such that $z$ and $z'$ are
opposite in the spherical residue $\proj_{R'_\eps}(R_\eps)$ (see
Section~\ref{sect:bldg:sph}). Set $x:=\op_\Sigma(z)$ and
$x':=\op_\Sigma(z')$. Notice that $x$ and $x'$ belong to
$R_{-\eps} \cap \Sigma$. We also define $y:= \proj_{R_\eps}(x)$.

We have $U(R_{-\eps}) \leq U(x')$ since $x' \in R_{-\eps}$.
Moreover, $U(R_{-\eps}) \cap \mathrm{Stab}_G(R_{\eps})$ fixes $y =
\proj_{R_\eps}(x)$ because this group fixes $x$ and stabilizes
$R_\eps$; hence, $U(R_{-\eps}) \cap \mathrm{Stab}_G(R_{\eps}) \leq
\mathrm{Stab}_G(y)$. Therefore, we have $U(R_{-\eps}) \cap
\mathrm{Stab}_G(R_{\eps}) \leq U(x') \cap \mathrm{Stab}_G(y) =
U(x', y)$, where the latter equality follows Lemma~\ref{U_w}.

We now choose a minimal gallery $y=y_0, \dots, y_j =
\proj_{R_\eps}(z'), \dots, y_n = z'$ joining $y$ to $z'$ (see
Section~\ref{sect:bldg:proj}). Hence, for all $0 \leq i \leq n$,
we have $y_i \in R_\eps$ if and only if $i \leq j$. (Notice that
$j = 0$ or $j=n$ is possible.)

For each $1 \leq i \leq n$, let $\beta_i$ be the twin root of
$\Sigma$ which contains $y_{i-1}$ but not $y_i$. Thus $\{\beta_1,
\dots, \beta_n\}$ is the set of twin roots containing $x'$ and $y$
or equivalently, $y$ but not $z'$ since $z' = \op_\Sigma(x')$. By
definition, we have $U(x', y) = \langle U_{\beta_i} | 1 \leq i
\leq n \rangle$. By Lemma \ref{U_w} this group coincides with the
product $U_{\beta_1}.U_{\beta_2} \dots U_{\beta_n}$.

Now we observe that by the definition of $y$, $y_j$ and $z'$ and
by Lemmas~\ref{Abr prop4} and~\ref{isom parallel}, we have
$$\proj_{R'_\eps}(y)=\proj_{R'_\eps}(y_j)=z'   \qquad \text{and}
\qquad \proj_{\proj_{R_\eps}(R_{-\eps})}(z')=
\proj_{\proj_{R_\eps}(R_{-\eps})}(y_j)=y.$$
In view of the
properties of projections (see Section~\ref{sect:bldg:proj}), this
implies
$$\begin{array}{rcl}
\{\beta \in \Phi^\Sigma(R_\eps) | \; \proj_{R_\eps}(R_{-\eps})
\subset \beta\} & = & \{\beta \in \Phi^\Sigma | \; y \in \beta
\text{ and } y_j \not \in \beta\}\\
&=& \{\beta_1, \dots, \beta_j\}\\
& = &  \{\beta \in \Phi^\Sigma | \; y \in \beta \text{ and }
\op_\Sigma(y_j) \in \beta\}
\end{array}$$
and
$$\begin{array}{rcl}
\{\beta \in \Phi^\Sigma | \; R_\eps \cap \Sigma \subset \beta
\text{ and } R_{-\eps} \cap \Sigma \subset \beta\} & = & \{\beta
\in \Phi^\Sigma | \; R_\eps \cap \Sigma \subset \beta \text{ and }
R'_\eps \cap \Sigma \subset \beta\}\\
& = & \{\beta \in \Phi^\Sigma | \; y_j \in \beta \text{ and } z'
\not \in \beta\}\\
& = & \{\beta_{j+1}, \dots, \beta_n\}\\
&=&\{\beta \in \Phi^\Sigma | \; y_j \in \beta \text{ and } x' \in
\beta\}.
\end{array}$$
We deduce from this, together with Lemma~\ref{U_w}, that
$$\widetilde{U}^\Sigma(R_\eps, R_{-\eps})=
\la U_i | \;  1 \leq i \leq j \ra =  U(y, \op_\Sigma(y_j))=U_1 \dots
U_j$$
and
$$U^\Sigma(R_\eps, R_{-\eps})= \la U_i | \;  j+1 \leq i \leq n \ra =  U(x', y_j)=U_{j+1} \dots
U_n.$$

In summary, we have shown that
$$\begin{array}{rcl}
U(R_{-\eps}) \cap \mathrm{Stab}_G(R_{\eps}) & \leq & U(x', y)\\
& = &\big(U_{\beta_1} \dots U_{\beta_j}\big).\big(U_{\beta_{j+1}}
\dots U_{\beta_n}\big)\\
& = & \widetilde{U}^\Sigma(R_\eps, R_{-\eps}).U^\Sigma(R_+,
R_-).
\end{array}$$
\end{proof}

The lemma implies that, if for $\eps = +$ or $-$ we have
$\proj_{R_\eps}(R_{-\eps})=R_\eps$ (in particular, if $R_+$ and
$R_-$ are parallel) and hence $\widetilde{U}^\Sigma(R_\eps,
R_{-\eps})=\{1\}$, then the group $U^\Sigma(R_+, R_-)$ is
independent of the twin apartment $\Sigma$. In that case, we may
omit the superscript $\Sigma$ and we shall write $U(R_+, R_-)$
rather than $U^\Sigma(R_+, R_-)$.

\begin{remark}\label{rem:Utilde} The group $\widetilde{U}^\Sigma(R_\eps, R_{-\eps})$
defined in the previous lemma actually coincides with \sloppy
$U(\proj_{R_\eps}(R_{-\eps})) \cap L^\Sigma(R_\eps)$. This can be
seen as follows. First notice that $(L^\Sigma, (U_\beta)_{\beta
\in \Phi^\Sigma(R_\eps)})$ is a twin root datum. Hence the group
$\mathrm{Stab}_{L^\Sigma(R_\eps)}(\proj_{R_\eps}(R_{-\eps}))$ has
a Levi decomposition in $L^\Sigma(R_\eps)$ by Proposition
\ref{Levi}. The above claim is easily deduced from that fact:
actually, the group $\widetilde{U}^\Sigma(R_\eps, R_{-\eps}) =
U(\proj_{R_\eps}(R_{-\eps})) \cap L^\Sigma(R_\eps)$ is nothing but
the unipotent radical of
$\mathrm{Stab}_{L^\Sigma(R_\eps)}(\proj_{R_\eps}(R_{-\eps}))$ with
respect to the above-mentioned twin root datum. We will not need
that fact here.
\end{remark}

The following is a consequence of the proof  of the previous
lemma.

\begin{cor}\label{nilpotence of U}
Let $R_+ \subseteq \Delta_+$ and $R_- \subseteq \Delta_-$ be
spherical residues. Then, for each $\eps \in \{+, -\}$ there exist
chambers $x_+ \in R_+$ and $x_- \in R_-$ such that $U(R_{-\eps})
\cap \mathrm{Stab}_G(R_\eps) = U(x_+, x_-)$. In particular, if all
root subgroups are nilpotent (i.e. if Condition (P1) holds), then
$U^\Sigma(R_+, R_-)$ is nilpotent, where $\Sigma$ is any twin
apartment intersecting $R_+$ and $R_-$.
\end{cor}
\begin{proof} The first statement was proved along the way. The second
statement is a consequence of the first, using also Axiom (TRD1)
and Lemma \ref{U_w}.
\end{proof}

\subsubsection{The intersection of finite type parabolics of opposite signs}
\label{sect:doubleLevi:result}

 We are now able to prove the Levi decomposition of intersections
 of finite type parabolic subgroups of opposite signs.

\begin{prop}\label{Levi for parallel}
Let $R_+ \subseteq \Delta_+$ and $R_- \subseteq \Delta_-$ be
spherical residues of $\Delta(G)$. Let $\Sigma$ be a twin
apartment containing $R_+$ and $R_-$. For each sign $\eps$ set
$R^\circ_\eps:= \proj_{R_\eps}(R_{-\eps})$. Then, for all $\eps,
\eps' \in \{+, -\}$ we have
$$\begin{array}{rcl}
\mathrm{Stab}_G(R_+) \cap \mathrm{Stab}_G(R_-) & = & L^\Sigma(R^\circ_{\eps'}) \ltimes U(R^\circ_+, R^\circ_-)\\
 & = & L^\Sigma(R^\circ_{\eps'}) \ltimes
 \big(\widetilde{U}^\Sigma(R_{-\eps}, R_\eps).
 U^\Sigma(R_+, R_-).\widetilde{U}^\Sigma(R_\eps, R_{-\eps})\big).\end{array}$$
\end{prop}
\begin{proof} Notice first that $\mathrm{Stab}_G(R_+) \cap
\mathrm{Stab}_G(R_-) = \mathrm{Stab}_G(R^\circ_+) \cap
\mathrm{Stab}_G(R^\circ_-)$. Moreover, since $R^\circ_+$ and
$R^\circ_-$ are parallel, we deduce from Lemma \ref{Abr prop4} and
Proposition \ref{parallelism}(ii) that $L^\Sigma(R^\circ_+) =
L^\Sigma(R^\circ_-)$.

Now the inclusion $L^\Sigma(R^\circ_{\eps'}) . U(R^\circ_+,
R^\circ_-) \leq \mathrm{Stab}_G(R_+) \cap \mathrm{Stab}_G(R_-)$ is
clear. On the other hand, we have the following:
$$\begin{array}{rcl}
\mathrm{Stab}_G(R^\circ_+) \cap \mathrm{Stab}_G(R^\circ_-) & = &
\big(L^\Sigma(R^\circ_\eps).U(R^\circ_\eps) \big) \cap
\mathrm{Stab}_G(R^\circ_{-\eps})\\
 & \leq & L^\Sigma(R^\circ_\eps). \big( U(R^\circ_\eps) \cap
\mathrm{Stab}_G(R^\circ_{-\eps}) \big)\\
 & = & L^\Sigma(R^\circ_{\eps'}). U(R^\circ_\eps ,
 R^\circ_{-\eps}),
\end{array}$$
where the last equality follows from Lemma \ref{U(R_+, R_-)}. This
proves that \sloppy $\mathrm{Stab}_G(R_+) \cap
\mathrm{Stab}_G(R_-) = L^\Sigma(R^\circ_{\eps'}). U(R^\circ_\eps,
R^\circ_{-\eps})$.

Now $L^\Sigma(R^\circ_{\eps'})=L^\Sigma(R^\circ_\eps)$ intersects
$U(R^\circ_\eps)$ trivially and normalizes that group by
Proposition \ref{Levi}. Since $L^\Sigma(R^\circ_{\eps'}) \leq
\mathrm{Stab}_G(R^\circ_{-\eps})$ we also see that
$L^\Sigma(R^\circ_{\eps'})$ normalizes
$\mathrm{Stab}_G(R^\circ_{-\eps})$. Therefore,
$L^\Sigma(R^\circ_{\eps'})$ normalizes $U(R^\circ_\eps) \cap
\mathrm{Stab}_G(R^\circ_{-\eps}) = U(R^\circ_\eps ,
R^\circ_{-\eps})$ and intersects the latter group trivially. This
proves the first equality of the lemma.

In order to establish the second equality, we first notice that
$\widetilde{U}^\Sigma(R_{-\eps}, R^\circ_\eps) =
\widetilde{U}^\Sigma(R_{-\eps}, R_\eps)$ by the definition of
these groups. Now, (an argument as in) the proof of the previous
lemma shows that
$$\begin{array}{rcl}
U(R^\circ_+, R^\circ_-) & = & \widetilde{U}^\Sigma(R_{-\eps},
R^\circ_\eps).U(R^\circ_\eps, R_{-\eps})\\
& = & \widetilde{U}^\Sigma(R_{-\eps},R_\eps) . U^\Sigma(R_\eps,
R_{-\eps}).\widetilde{U}^\Sigma(R_\eps, R_{-\eps}),\end{array}$$
from which the conclusion follows.
\end{proof}

\begin{remark} The preceding proposition is proved in \cite{Re99},
Section 6.3.4, under an additional assumption called (NILP),
defined in \emph{op. cit.}, Section 6.3. Our proof shows that this
extra assumption is not necessary for the result to hold.
\end{remark}

\begin{cor}\label{G_f <=> finite}
Let $(W, S)$ be a Coxeter system such that $S$ is finite. Let
\sloppy $\mathcal{D}=(G, (U_\phi)_{\phi \in \Phi(W, S)})$ be a
twin root datum such that each $U_\phi$ is finite and such that $H
:= \bigcap_{\phi \in \Phi(W,S)} N_G(U_\phi)$ is finite. Then the
set of all bounded subgroups coincides with the set of all finite
subgroups of $G$.
\end{cor}
\begin{proof} Let $\Delta$ be the twin building associated with
$\mathcal{D}$. The fact that each finite subgroup of $G$ is
bounded is an immediate consequence of Proposition \ref{fixed
point}. In order to prove that a bounded subgroup is finite, it
suffices to prove that given a pair $R_+$, $R_-$ of spherical
residues of opposite signs, then the group $\mathrm{Stab}_G(R_+)
\cap \mathrm{Stab}_G(R_-)$ is finite. Our hypotheses imply that
every Levi factor of spherical type is finite. Hence, by
Proposition \ref{Levi for parallel}, it just remains to show that
$U(\proj_{R_+}(R_-), \proj_{R_-}(R_+))$ is finite. But this
follows again from our hypotheses in view of Lemma \ref{U_w} and
Corollary \ref{nilpotence of U}.
\end{proof}

\section{Maximal bounded subgroups}\label{max sub sect}

\subsection{The main characterization}

Let $\mathcal{D}=(G, (U_{\alpha})_{\alpha \in \Phi})$ be a twin
root datum of type $(W,S)$. By a maximal bounded subgroup of $G$,
we mean a bounded subgroup which is not properly contained in any
other bounded subgroup of $G$. Let $M$ be such a bounded subgroup.
By definition $M$ is the intersection of two finite type parabolic
subgroups of opposite signs. The following theorem shows that
there exists two canonical finite type parabolic subgroups $P^M_+$
and $P^M_-$ such that $M= P^M_+ \cap P^M_-$. Case~(i) of the
theorem corresponds to the case where $P^M_+$ and $P^M_-$ are
opposite; the group $M$ is then the common Levi factor of $P^M_+$
and $P^M_-$. To some extent, this is the generic case (see
Proposition~\ref{max sub cor} and Remark~\ref{rem:(i)generic}
below).

Before stating the theorem, we need one more notation. Let $\tb$
be the twin building associated with $\mathcal{D}$. Given a
subgroup $M \leq G$, we denote by $\mathcal{S}_\eps(M)$ the set of
all spherical residues of $\Delta_\eps$ stabilized by $M$, where
$\eps \in \{+,-\}$.

\begin{theo}\label{max sub}
Let \TRD{} be a twin root datum of type (W,S) and let $M \leq G$
be a maximal bounded subgroup. Then one of the following
holds~:\\
\begin{tabularx}{\linewidth}{rX}
(i) & for $\eps \in \{+, -\}$ the set $\mathcal{S}_\eps(M)$
consists of a unique element $R_\eps$, which is a maximal
spherical residue of $\Delta_\eps$; moreover $R_+$ and $R_-$ are
opposite in
$\Delta(G)$;\\
(ii) & for $\eps \in \{+, -\}$ the set $\mathcal{S}_\eps(M)$
possesses two distinguished elements $R_\eps$ and
$\overline{R}_\eps$ such that for every $T_\eps \in
\mathcal{S}_\eps(M)$ we have $R_\eps \subseteq T_\eps \subseteq
\overline{R}_\eps$ and $\proj_{T_\eps}(T_{-\eps})=R_\eps$ and
$T_+$ and $T_-$ are not opposite; moreover $\overline{R}_\eps$ is
the only maximal spherical residue containing $R_\eps$.
\end{tabularx}
In both cases we have $M = \mathrm{Stab}_G(R_+) \cap
\mathrm{Stab}_G(R_-)$.

Conversely, let $R_+ \subseteq \Delta_+$ and $R_- \subseteq
\Delta_-$ be spherical residues such that either of the following conditions holds~:\\
\begin{tabularx}{\linewidth}{rX}
(i$'$) & $R_+$ and $R_-$ are maximal spherical and opposite;\\
(ii$'$) & $R_+$ and $R_-$ are parallel; moreover, for each $\eps
\in \{+, -\}$ the residue $R_\eps$ is properly contained in a
unique maximal spherical residue $\overline{R}_\eps$ and we have
$proj_{\overline{R}_\eps}(\overline{R}_{-\eps}) = R_\eps$.
\end{tabularx}
Then $M := \mathrm{Stab}_G(R_+) \cap \mathrm{Stab}_G(R_-)$ is a
maximal bounded subgroup.
\end{theo}
\begin{proof} Let $M \leq G$ be a maximal bounded subgroup. For
$\eps \in \{+, -\}$ let $R_\eps \in \mathcal{S}_\eps (M)$.

Assume first that $R_+$ and $R_-$ are opposite. Hence, $M =
\mathrm{Stab}_G(R_+) \cap \mathrm{Stab}_G(R_-)$ which implies by
Proposition~\ref{Levi} that $M$ does not stabilize any residue
properly contained in $R_+$ or $R_-$. Moreover, the maximality of
$M$ implies that $R_+$ and $R_-$ are maximal spherical residues.
Let now $T_\eps \in \mathcal{S}_\eps(M)$. Then $M$ stabilizes
$\proj_{R_+}(T_\eps)$ and $\proj_{R_-}(T_\eps)$. Since these
cannot be properly contained in $R_+$ and $R_-$ respectively, we
conclude from Corollary~\ref{proj max sph} that $T_\eps = R_\eps$.
Hence, we are in Case (i). Notice that our discussion also proves
the converse statement in the case (i$'$).

We now assume that $R_+$ and $R_-$ are not opposite. For $\eps \in
\{+, -\}$ we have $\proj_{R_\eps}(R_{-\eps}) \in
\mathcal{S}_\eps(M)$. Moreover, we know by Lemma \ref{iff
opposite} that $\proj_{R_+}(R_-)$ and $\proj_{R_-}(R_+)$ are not
opposite (this follows also from the first part of the present
proof). Therefore, up to replacing $R_\eps$ by
$\proj_{R_\eps}(R_{-\eps})$ we may assume that $R_+$ and $R_-$ are
parallel. Let $J_+$ and $J_-$ be their respective types. Since
$R_+$ and $R_-$ are not opposite, Lemma \ref{Abr prop4} and
Proposition \ref{parallelism} show that $J_\eps$ is not a maximal
spherical subset of $S$.

Let now $s \in S \backslash J_\eps$ such that $J_\eps \cup \{s\}$
is spherical and let $R^s_\eps$ be the residue of type $J_\eps
\cup \{s\}$ containing $R_\eps$. Set also $T :=
\proj_{R_\eps^s}(R_{-\eps})$. We now prove that $T = R_\eps$.

Assume on the contrary that $T \neq R_\eps$. Then Lemma
\ref{compose proj} and Proposition \ref{parallelism} imply that
$R_\eps$ and $T$ are opposite in $R_\eps^s$. Let $\Sigma$ be a
twin apartment containing $R_+$ and $R_-$. There exists a twin
root $\alpha = (\alpha_+, \alpha_-)$ of $\Sigma$ such that $T
\subseteq \alpha_\eps$, $R_\eps \subseteq -\alpha_\eps$ (hence the
reflection $s_\alpha$ of $\Sigma$ stabilizes $R_\eps^s$) and
$R_{-\eps} \subseteq \alpha_{-\eps}$. Now Proposition \ref{Levi}
implies that $U_\alpha$ acts freely on the residues opposite $T$
in $R_\eps^s$. In particular, we have $U_\alpha \cap M = \{1\}$
since $M$ stabilizes $R_\eps$. Therefore, the group $\langle
U_\alpha \cup M \rangle$ contains $M$ properly, and it stabilizes
$R_\eps^s$ and $R_{-\eps}$ by construction. In other words, we
have $M \lneq \langle U_\alpha \cup M \rangle$ is a bounded
subgroup, which contradicts the maximality of $M$. This proves
that $T =R_\eps$ as claimed.

Let $S_\eps$ be the set of all $s \in S\backslash J_\eps$ such
that $J_\eps \cup \{s\}$ is spherical. Let $\overline{R}_\eps$ be
the residue of type $J_\eps \cup S_\eps$ containing $R_\eps$. We
now prove that $\overline{R}_\eps$ is spherical.

To this end we consider, as before, a twin apartment $\Sigma$
containing $R_+$ and $R_-$. Let $R'_\eps := \op_\Sigma (R_\eps)$.
Choose a chamber $x \in R_\eps$ and a chamber $z$ which is
opposite to $\proj_{R'_\eps}(x)$ in $R'_\eps$. Set
$y:=\proj_{\overline{R}_\eps}(z)$. Our aim is to apply the
criterion of sphericity of Lemma~\ref{criterion of spher} to $x$
and $y$. Hence, let $j \in J_\eps \cup S_\eps$ and denote by
$\pi_j$ the $j$-panel containing $x$.

We know from Lemma \ref{Abr prop4} that $R_\eps$ and $R'_\eps$ are
parallel. By Lemma \ref{compose proj}, this implies that $R_\eps$
and $\proj_{\overline{R}_\eps}(R'_\eps)$ are parallel. By
Lemma~\ref{isom parallel}, we see that $x$ and $\proj_{R_\eps}(y)$
are opposite in $R_\eps$ and we conclude from Lemma~\ref{criterion
of spher} that if $j \in J_\eps$ then $x \neq
\proj_{\pi_j}(y)=\proj_{\pi_j} (\proj_{R_\eps}(y))$.

Let us now assume that $j \in S_\eps$. We have already proved that
$\proj_{R^j_\eps}(R_{-\eps})=R_\eps$, which implies that
$\proj_{R^j_\eps}(R'_\eps)$ is opposite $R_\eps$ in $R_\eps^j$ by
Lemma \ref{Abr prop4}. Therefore, $x$ and $\proj_{R^j_\eps}(y)$
are opposite in $R_\eps^j$ which implies by Lemma \ref{criterion
of spher} that $x \neq \proj_{\pi_j}(y)=\proj_{\pi_j}
(\proj_{R^j_\eps}(y))$.

Finally, Lemma \ref{criterion of spher} applied to $x$ and $y$
insures that $\overline{R}_\eps$ is spherical, i.e. that $J_\eps
\cup S_\eps$ is a spherical subset of $S$. Moreover, it is clear
by the definition of $S_\eps$ that $J_\eps \cup S_\eps$ is
actually a \emph{maximal} spherical subset of $S$.

By maximality of $M$, we have  $M = \mathrm{Stab}_G(R_+) \cap
\mathrm{Stab}_G(R_-)$. Therefore, we see by Proposition \ref{Levi
for parallel} that $M$ does not stabilize any proper residue of
$R_\eps$ for $\eps \in \{+, -\}$. Moreover, the same result
implies that if $R$ is a residue contained in $\overline{R}_\eps$,
then $R$ is stabilized by $M$ if and only if $R$ contains
$R_\eps$.

Let now $T_\eps \in \mathcal{S}_\eps(M)$. Then
$\proj_{R_\eps}(T_\eps) \in \mathcal{S}_\eps(M)$ and so
$\proj_{R_\eps}(T_\eps)=R_\eps$. Therefore, $R_\eps$ and
$\proj_{T_\eps}(R_\eps)$ are parallel. By the previous paragraph
together with Proposition \ref{parallelism}, this implies that
$R_\eps = \proj_{T_\eps}(R_\eps)$, or in other words that $R_\eps
\subseteq T \subseteq \overline{R}_\eps$.

Going now back to the first argument in the proof of (ii) above,
we conclude moreover that $\proj_{T_\eps}(T_{-\eps})=R_\eps$. The
fact that $T_+$ and $T_-$ are not opposite is now obvious.

It remains to prove the converse statement. For (ii$'$), let
$R_+$, $R_-$, $\overline{R}_+$ and $\overline{R}_-$ be as in
(ii$'$) and define $M := \mathrm{Stab}_G(R_+) \cap
\mathrm{Stab}_G(R_-)$. Mimicking some of the arguments above we
can prove again that if $R$ is a residue contained in
$\overline{R}_\eps$, then $R$ is stabilized by $M$ if and only if
$R$ contains $R_\eps$. From this, we deduce as above that $R_\eps
\subseteq T \subseteq \overline{R}_\eps$ whenever $T \in
\mathcal{S}_\eps(M)$. Therefore, if $M \leq M_1$ and $M_1$ is
bounded, then there exists $T_\eps \in \mathcal{S}_\eps(M_1)$ such
that $R_\eps \subseteq T_\eps \subseteq \overline{R}_\eps$. But
our hypotheses then imply that $\proj_{T_\eps}(T_{-\eps})=R_\eps$.
Therefore, we have $R_\eps \in \mathcal{S}_\eps(M_1)$ namely $M_1
\leq \mathrm{Stab}_G(R_\eps)$. Hence $M = M_1$. The proof is
complete.
\end{proof}

\subsection{Two specializations}

\subsubsection{Obstructions for Case~(ii) of Theorem~\ref{max sub}}
In many interesting situations, only Case~(i) of Theorem~\ref{max
sub} occurs. The next result gives sufficient conditions on the
Coxeter system $(W,S)$ which imply that Case~(ii) never happens.

\begin{prop}\label{max sub cor}
Assume that the Coxeter system $(W,S)$ satisfies one of the
following conditions~:\\
\begin{tabularx}{\linewidth}{lX}
(R1) & for all $s, t \in S$, the order of $st$ is not equal to
3;\\
(R2) & for any pair $J, K$ of spherical subsets of $S$ such that
$J$ is properly contained in $K$ and $K$ is maximal spherical,
there exists an $s \in S \backslash K$ such that $J \cup \{s\}$ is
spherical but $K \cup
\{s\}$ is not;\\
(R3) & for every $j \in S$, there exists a unique maximal
spherical subset $J$ of $S$ such that $j \in J$.
\end{tabularx}\\
Then the case (ii) does not occur in the previous theorem.
\end{prop}
\begin{proof} We keep the notation of the proof of Theorem \ref{max sub}. It
is clear that if (ii) holds in that theorem, then (R2) fails by
choosing $J = J_\eps$ and $K=J_\eps \cup S_\eps$ where $\eps = +$
or $-$.

Now let $\Sigma$ be a twin apartment containing $R_+$ and $R_-$ as
in the proof above, let $R'_\eps:=\op_\Sigma(R_{-\eps})$ and let
$J''_\eps$ be the type of
$R''_\eps:=proj_{\overline{R}_\eps}(R'_\eps)$. We know by
Lemma~\ref{Abr prop4} that $R''_\eps$ is opposite $R_\eps$ in
$\overline{R}_\eps$ and, moreover, that $\overline{R}_\eps$ is the
\emph{unique} maximal spherical residue containing $R_\eps$. On
the other hand, since $\overline{R}_+$ and $\overline{R}_-$ are
not opposite, we have $R'_\eps \neq R''_\eps$, and these two
distinct residues are parallel. Therefore, it follows from
Proposition \ref{parallelism} that there exists $s \in S
\backslash(J_\eps \cup S_\eps)$ such that $J''_\eps \cup \{s\}$ is
spherical. In particular, $J_\eps \cup S_\eps$ is \emph{not} the
unique maximal spherical subset of $S$ containing $J''_\eps$.
Hence, (R3) fails and furthermore, we have $J_\eps \neq J''_\eps$.
Since $J_\eps$ and $J''_\eps$ are the respective types of two
opposite residues of $\overline{R}_\eps$, the latter inequality
also implies that (R1) fails, using \cite{Sc95}, Proposition~5.2.3
and the fact that there are no Moufang $n$-gons for odd $n$
greater than 3.
\end{proof}

\begin{remark}\label{rem:(i)generic}
\begin{enumerate}

\item Condition (R2) in the previous corollary is also equivalent
to the following~:\\
\begin{tabularx}{\linewidth}{rX}
(R2$'$) & for every non-maximal spherical subset $J$ of $S$ there
exists (at least) two distinct maximal spherical subsets $K_1$ and
$K_2$ of $S$ containing $J$.
\end{tabularx}

\item All affine and compact hyperbolic Coxeter diagrams satisfy
Condition (R2) (notice that (R2) is empty for $\tilde{A}_1$~ and
so obviously satisfied; actually $\tilde{A}_1$ also satisfies (R1)
and (R3)).

\item Condition (R3) in the previous corollary is also
equivalent to each of the following ones~:\\
\begin{tabularx}{\linewidth}{rX}
(R3$'$) & for every maximal spherical subset $J$ of $S$ and all
pairs $j, s$ with $j \in J$ and $s \in S \backslash J$, the order
of $sj$ is infinite;\\
(R3$''$) & there is a partition $S = S_1 \cup \dots \cup S_n$ of
$S$ into spherical subsets such that the order of $st$ is infinite
whenever $s \in S_i$, $t \in S_j$ and $i \neq j$.
\end{tabularx}

\end{enumerate}
\end{remark}

\subsubsection{A group theoretic description of Case~(ii)}
We end this section with a lemma which is will be used in
the proof of Theorem \ref{main}.

\begin{lem}\label{char type (ii)}
Let \TRD{} be a twin root datum of type (W,S) which satisfies
Conditions (P1) and (P2) of Section \ref{sect:TRD:condP}. Let $M
\leq G$ be a maximal bounded subgroup. A necessary and sufficient
condition for $M$ to have type (ii) in Theorem \ref{max sub}, is
that
$$M = U \rtimes (M \cap M'),$$
where $U$ is a nontrivial nilpotent group and $M'$ is a maximal
bounded subgroup different from $M$.
\end{lem}
\begin{proof} We first prove that the condition is necessary. We keep the
notation of Theorem~\ref{max sub} and assume that $M$ satisfies
Condition (ii). Let $\Sigma$ denote a twin apartment containing
$R_+$ and $R_-$. By Proposition \ref{Levi for parallel}, we have
$$\begin{array}{rcl}
M & = & \mathrm{Stab}_G(R_+) \cap \mathrm{Stab}_G(R_-)\\
  & = & \mathrm{Stab}_G(\overline{R}_+) \cap \mathrm{Stab}_G(R_-)\\
  & = & L^\Sigma(R_+).\widetilde{U}^\Sigma(\overline{R}_+,
  R_-).U(\overline{R}_+, R_-). \widetilde{U}^\Sigma(R_-,
  \overline{R}_+)\\
  & = & L^\Sigma(R_+).\widetilde{U}^\Sigma(\overline{R}_+,
  R_-).U(\overline{R}_+, R_-),\end{array}$$
where the last equality follows from
$\proj_{R_-}(\overline{R}_+)=R_-$ which implies
$\widetilde{U}^\Sigma(R_-, \overline{R}_+)=\{1\}$.

Now, set $U:=U(\overline{R}_+, R_-)$ and
$M':=L^\Sigma(\overline{R}_+)$. By Theorem \ref{max sub}, the
group $M'$ is a maximal bounded subgroup.  Clearly, the group $U$
is nilpotent by \ref{nilpotence of U} and nontrivial since
$\op_\Sigma(R_-)\cap \overline{R}_+ = \emptyset$ (see the proof of
Proposition \ref{max sub cor}). Moreover, we have $M' \cap U
=\{1\}$ since $U \leq U(\overline{R}_+)$. On the other hand, we
also have $L^\Sigma(R_+).\widetilde{U}^\Sigma(\overline{R}_+, R_-)
\leq M'$ by definition. Therefore, we have
$L^\Sigma(R_+).\widetilde{U}^\Sigma(\overline{R}_+, R_-) = M \cap
M'$ and hence, $M = (M \cap M').U$.

It remains to prove that $M\cap M'$ normalizes $U$. Using again
the fact that $\widetilde{U}^\Sigma(R_-, \overline{R}_+)$ is
trivial, we deduce from Lemma \ref{U(R_+, R_-)} that $U =
U(\overline{R}_+) \cap \mathrm{Stab}_G(R_-)$. Now the desired
conclusion follows since $M$ normalizes $\mathrm{Stab}_G(R_-)$
(because $M \leq \mathrm{Stab}_G(R_-)$) and $M'$ normalizes
$U(\overline{R}_+)$ (by Proposition \ref{Levi}).

We now show that the condition is sufficient. Assume that $M$ has
type (i) in Theorem~\ref{max sub} and that the condition of the
lemma is satisfied. Let $R_+$ (resp. $R_-$) be the unique element
of $\mathcal{S}_+(M)$ (resp. $\mathcal{S}_-(M)$) and let $T_+ \in
\mathcal{S}_+$. Since $M$ and $M'$ are distinct (or since $U$ is
nontrivial), the residues $R_+$ and $T_+$ are distinct. Now, $M
\cap M'$ stabilizes $R:=\proj_{R_+}(T_+)$, which is properly
contained in $R_+$ by Corollary \ref{proj max sph}. In particular,
the group $U$ does not act trivially on $R_+$ since $M = (M \cap
M').U$ does not stabilize $R$.

As before, let $\Sigma$ be a twin apartment intersecting $R_+$ and
$R_-$. Then $M = L^\Sigma(R_+)$ and we know that $(M,
(U_\alpha)_{\alpha \in \Phi^\Sigma(R_+)})$ is a twin root datum of
spherical type $(W_J, J)$ for some $J \subseteq S$. Since $U$ is
normal in $M$ and since it does not act trivially on $R_+$, we
deduce from \cite{Bo81}, Theorem 5 of Chapter IV that there exists
a residue $R' \subseteq R_+$ such that $U$ contains the group
$M_1$ generated by all $U_\beta$ with $\beta \in \Phi^\Sigma(R')$.
But $\beta \in \Phi^\Sigma(R')$ implies $-\beta \in
\Phi^\Sigma(R')$ and therefore $M_1$ is generated by subgroups of
the form $\langle U_\beta \cup U_{-\beta} \rangle$. Since each
group of the latter form is perfect by hypothesis, the group $M_1$
itself is perfect which contradicts the fact that it is contained
in the nilpotent group $U$. This concludes the proof of the lemma.
\end{proof}

\section{The reduction theorem for isomorphisms which preserve bounded subgroups}
In this section we state and prove a general theorem concerning
isomorphisms between groups endowed with twin root data, which
preserve bounded subgroups. Roughly speaking, it says that, under
certain conditions, the isomorphism problem for groups with twin
root data reduces to the isomorphism problem for groups of finite
type. The main results of this paper will be deduced from it in
the following two sections.

\subsection{$\mathfrak{E}$-rigidity of twin root data}

In order to make the statement of this theorem as precise and
concise as possible, we introduce some additional terminology.

Let $\mathfrak{E}$ be a collection of twin root data. A twin root
datum $\mathcal{D}$ is called
\textbf{$\mathfrak{E}$-rigid} if the following holds (see
Section~\ref{sect:TRD:defPrinc} for the definition of $G^\mathcal{D}$):\\
\begin{center}\vspace{-.5cm}
If $\mathcal{D}' \in \mathfrak{E}$, then any isomorphism of
$G^\mathcal{D}$ to $G^{\mathcal{D}'}$ induces an isomorphism of
$\mathcal{D}$ to $\mathcal{D}'$.
\end{center}

Let $\mathcal{D}=(G, (U_{\alpha})_{\alpha \in \Phi(W,S)})$ be a
twin root datum, let $\Delta$ be the associated twin building, let
$\Sigma$ be the fundamental twin apartment and let $c$ be a
chamber of $\Sigma$. For each subset $J$ of $S$, we set $L_J :=
L^\Sigma(\mathrm{Res}_J(c))$ and $\mathcal{D}_J:=(L_J,
(U_{\alpha})_{\alpha \in \Phi^\Sigma(\mathrm{Res}_J(c))})$.

Given a collection $\mathfrak{E}$ as above, then we denote by
$\mathfrak{E}_\mathrm{sph}$ the collection of all twin root data
of the form $\mathcal{D}_J$ such that $\mathcal{D} \in
\mathfrak{E}$ is of type $(W^\mathcal{D}, S^\mathcal{D})$ and $J$
is a maximal spherical subset of $S^\mathcal{D}$. A twin root
datum $\mathcal{D} \in \mathfrak{E}$ of type $(W^\mathcal{D},
S^\mathcal{D})$ is called \textbf{$\mathfrak{E}$-locally rigid} if
for every every maximal spherical subset $J$ of $S^\mathcal{D}$,
the twin root datum $\mathcal{D}_J$ is
$\mathfrak{E}_{\mathrm{sph}}$-rigid.

\subsection{The result}

\begin{theo}\label{main}
Let $\mathcal{D}$ and $\mathcal{D}'$ be two twin root data which
satisfy Conditions (P1), (P2) and (P3), and whose types are
Coxeter systems of finite rank. Assume that $\mathcal{D}$ is
$\{\mathcal{D}'\}$-locally rigid. If $\xi : G^\mathcal{D} \to
G^{\mathcal{D}'}$ is an isomorphism which maps bounded subgroups
of $G^\mathcal{D}$ to bounded subgroups of $G^{\mathcal{D}'}$,
then $\xi$ induces an isomorphism of $\mathcal{D}$ to
$\mathcal{D}'$.
\end{theo}
\begin{proof} Let $M$ be a maximal bounded subgroup of $G^\mathcal{D}$.
Then $\xi(M)$ is a maximal bounded subgroup of ${G}^\mathcal{D'}$
by the hypothesis on $\xi$. Moreover, it follows from Lemma
\ref{char type (ii)} that $M$ and $\xi(M)$ have the same `type',
where the `type' of $M$, resp. $\xi(M)$ is given by either Case
(i) or (ii) in Theorem \ref{max sub}.

Let now $\Sigma$ be a twin apartment of the twin building $\Delta$
associated with $\mathcal{D}$ and let $\alpha$ be a twin root of
$\Sigma$. Let also $\Delta'$ be the twin building associated with
$\mathcal{D'}$. Choose a maximal residue of spherical type $R$
intersecting $\Sigma$ and such that $\alpha \in \Phi^\Sigma(R)$.
By what we have just seen, we know that $\xi(L^\Sigma(R))$ has the
form $L^{\Sigma'}(R')$ for some apartment $\Sigma'$ of $\Delta'$
and some maximal residue of spherical type $R'$ which intersects
$\Sigma'$.

Since $\mathcal{D}$ is locally rigid, the restriction of $\xi$ to
$M = L^\Sigma(R)$ induces an isomorphism from the twin root datum
$(L^\Sigma(R), (U_\beta)_{\beta \in \Phi^\Sigma(R)})$ to the twin
root datum $(L^{\Sigma'}(R'), (U'_\beta)_{\beta \in
\Phi^{\Sigma'}(R')})$, where we have used superscript~$'$ to
denote root groups of $\mathcal{D}'$. We may assume without loss
of generality that
$$\{\xi(U_\beta) | \beta \in \Phi^\Sigma(R)\} = \{
U'_\beta | \beta \in \Phi^{\Sigma'}(R')\}.$$

Let $H = \Fix_{G^\mathcal{D}}(\Sigma)$ and let $H':=\xi(H)$. Since
\sloppy $H = \bigcap_{\beta \in \Phi^\Sigma(R)}
N_{L^\Sigma(R)}(U_\beta)$ (see Section~\ref{sect:TRD:parab}) we
deduce  $H '= \bigcap_{\beta \in \Phi^{\Sigma'}(R')}
N_{L^{\Sigma'}(R')}(U'_\beta)$. This implies that $H'$ is the
chamberwise stabilizer of $\Sigma'$ in $G(\mathcal{D'})$.

Let now $\gamma$ be a twin root of $\Sigma$ which does not belong
to $\Phi^\Sigma(R)$. Arguing as for $\alpha$, we obtain that
$\xi(U_\gamma) = U'_{\gamma'}$ where $\gamma'$ is a twin root of
$\Delta'$ which is contained in a twin apartment $\Sigma''$ whose
chamberwise stabilizer is $H'$. On the other hand, our hypotheses
imply that $H'$ fixes a unique twin apartment chamberwise (see
Lemma \ref{H not trivial}). In summary, we have shown that
$$\{\xi(U_\beta) | \beta \text{ is a twin root of $\Sigma$}\} = \{
U'_\beta | \beta \text{ is a twin root of $\Sigma'$}\}.$$

Now the conclusion follows from Proposition~\ref{thm:KM1}.
\end{proof}

\section{Kac-Moody groups over arbitrary fields}
In this section and in the following one, we apply Theorem
\ref{main} to the case of Kac-Moody groups over fields in the
strict sense.

It is known that a Kac-Moody group $G$ over a field $\K$ naturally
yields a twin root datum $\mathcal{D}=(G, (U_\alpha)_{\alpha \in
\Phi})$ which is locally split over $\K$ (namely which is locally
split over $(\K_\alpha)_{\alpha \in \Phi}$, where $\K_\alpha = \K$
for each $\alpha \in \Phi$). We have also mentioned that if $\K$
has cardinality at least 4, then the conditions (P1), (P2) and
(P3) of Section~\ref{sect:TRD:condP} are satisfied. In order to
apply Theorem \ref{main} to $G$, it remains to discuss the local
rigidity of the twin root datum $\mathcal{D}$. This is done by
using the classical theorems on isomorphisms between Chevalley
groups, but the arguments are slightly different according as the
ground field is finite or infinite.

\subsection{Finite fields vs. infinite fields}

The following result gives a handy criterion which distinguishes
between these two cases.

\begin{prop}\label{fi field <=> fi gen}
Let $G$ be a Kac-Moody group over a field $\K$. Then $G$ is
finitely generated if and only if $\K$ is finite.
\end{prop}
\begin{proof}
Let $\mathcal{D}=(G,(U_\alpha)_{\alpha \in \Phi})$ be the twin
root datum which is naturally associated with $G$. Let $H :=
\bigcap_{\alpha \in \Phi} N_G(U_\alpha)$ and let $\Pi \subset
\Phi$ be the (finite) set of simple roots in $\Phi=\Phi(W,S)$. It
is known (and easy to see) that $G$ is generated by the set
$S_{\mathcal{D}} := H \cup \bigcup_{\alpha \in \Pi} U_\alpha$.
Moreover, each $U_\alpha$ is isomorphic to the additive group $\K$
and the group $H$ is a `split $\K$-torus', namely it is isomorphic
to a direct product of finitely many copies of the multiplicative
group $\K^\times$.

If $\K$ is finite, then $S_\mathcal{D}$ is finite, whence $G$ is
finitely generated.

If $G$ is finitely generated, then $G$ is generated by a subset of
$S_\mathcal{D}$. By \cite{Ti87}, the defining relations satisfied
by the elements of $S_\mathcal{D}$ in the group $G$ involve only
the ring structure of $\K$. Since no infinite field is a finitely
generated ring, we deduce that $\K$ has to be finite.
\end{proof}

\subsection{The characteristic in the case of a finite ground field}

The following result will spare us to worry about the exceptional
isomorphisms between finite Chevalley groups.

\begin{prop}\label{characteristic}
Let $G$ be an infinite Kac-Moody group over a finite field $\K$ of
characteristic $p$. Let $q$ be a prime. Then $p=q$ if and only if
the set of orders of finite $q$-subgroups of $G$ has no finite
upper bound.
\end{prop}
\begin{proof}
Let $\mathcal{D}=(G,(U_\alpha)_{\alpha \in \Phi})$ be the twin
root datum which is naturally associated with $G$.

Assume first that $p=q$. We must show that $G$ possesses finite
$p$-subgroups of arbitrary large orders. Since $\K$ is finite, the
group $U_\alpha$, which is isomorphic to the additive group $\K$,
is finite for every $\alpha \in \Phi$. Since $G$ is infinite by
assumption, we deduce that the Coxeter group $W$ is infinite. Let
$\Delta = (\Delta_+, \Delta_-, \cod)$ be the twin building
associated with $\mathcal{D}$ and let $\Sigma$ be the fundamental
twin apartment. Since $\Delta$ is of non-spherical type, we can
find chambers $x_+$ and $x_-$ of $\Sigma$ such that $n:=
\ell(\cod(x_+, x_-))$ is arbitrarily large. On the other hand, we
know by Lemma \ref{U_w} that the group $U(x_+, x_-)$ may be
written as a product of the form $U_{\beta_1}.U_{\beta_2}\dots
U_{\beta_n}$ for certain twin roots $\beta_1, \dots, \beta_n$ of
$\Sigma$. Since $U_{\beta_i}$ is a finite $p$-group for each $1
\leq i \leq n$, it follows that $U(x_+, x_-)$ is a finite
$p$-group of order at least $p^n$ which yields the desired result.

We now assume that $p \neq q$. We must show that there is an
upper-bound on the possible orders of finite $q$-subgroups of $G$.
Let $Q \leq G$ be such a finite $q$-group. By Proposition
\ref{fixed point}, the group $Q$ is a bounded subgroup. Let $R_+
\subset \Delta_+$ and $R_- \subset \Delta_-$ be spherical residues
which are stabilized by $Q$. Up to replacing $R_+$ and $R_-$ by
$\proj_{R_+}(R_-)$ and $\proj_{R_-}(R_+)$ respectively, we may
assume that $R_+$ and $R_-$ are parallel. Let $\Sigma_Q$ be a twin
apartment containing $R_+$ and $R_-$. By Proposition \ref{Levi for
parallel}, we have $\mathrm{Stab}_G(R_+) \cap \mathrm{Stab}_G(R_-)
= L^{\Sigma_Q}(R_+) \ltimes U(R_+, R_-)$. Hence there is a
homomorphism $f : \mathrm{Stab}_G(R_+) \cap \mathrm{Stab}_G(R_-)
\to L^{\Sigma_Q}(R_+)$. On the other hand, the order of every
element of $U(R_+, R_-)$ is a power of $p$ by Lemma \ref{U_w} and
Corollary \ref{nilpotence of U}. Since $p \neq q$ we deduce that
$Q$ and $f(Q)$ are isomorphic and, hence, we may assume that $Q$
is contained in $L^{\Sigma_Q}(R_+)$. Now, up to conjugation by an
element of $G$, we may assume that $L^{\Sigma_Q}(R_+) = L_J$ for
some spherical subset $J$ of $S$. Thus $|Q| \leq \mathrm{max}
\{|L_J| \ | J \subseteq S \text{ spherical}\}$ (note that each
$L_J$ is a finite Chevalley group over $\K$). The desired
conclusion follows from the finiteness of $S$.
\end{proof}

\subsection{Isomorphisms of Kac-Moody groups}

We are now able to apply Theorem \ref{main} to Kac-Moody groups
over fields. The following theorem is the main result of the
introduction.

\begin{theo}\label{main over arbitrary fields}
Let $G$ and $G'$ be infinite Kac-Moody groups over fields $\K$ and
$\K'$ respectively, both of cardinality at least 4, and let
$\mathcal{D}$ and $\mathcal{D'}$ be the corresponding twin root
data. Let $\xi : G \to G'$ be a group isomorphism which maps
bounded subgroups of $G$ to bounded subgroups of $G'$. Then $\xi$
induces an isomorphism of $\mathcal{D}$ to $\mathcal{D}'$.
\end{theo}
\begin{proof}
We have to show that $\xi$ induces an isomorphism of $\mathcal{D}$
to $\mathcal{D}'$. By Theorem \ref{main}, it suffices to show that
$\mathcal{D}$ is $\{\mathcal{D}'\}$-locally rigid. Let $(W, S)$
and $(W', S')$ be the respective types of $\mathcal{D}$ and
$\mathcal{D}'$, let $J$ and $J'$ be maximal spherical subsets of
$S$ and $S'$ and let $\varphi : L_J \to L_{J'}$ be a group
isomorphism. We have to show that $\varphi$ induces an isomorphism
of $\mathcal{D}_J$ to $\mathcal{D}'_{J'}$.

By the definition of a Kac-Moody group, we know that $L_J$ and
$L_{J'}$ are Chevalley groups over $\K$ and $\K'$ respectively. Up
to replacing $L_J$ (resp. $L_{J'}$) by its derived subgroup modulo
its center and $\mathcal{D}_J$ (resp. $\mathcal{D}'_{J'}$) by its
\emph{reduction} (see \cite{Ti92}, Section 3.3 and \cite{CM03b},
Section 3.13), we may assume that $L_J$ (resp. $L_{J'}$) is an
adjoint Chevalley group.

Since $G$ and $G'$ are isomorphic, it follows from Proposition
\ref{fi field <=> fi gen} that $\K$ and $\K'$ are either both
finite or both infinite. Moreover, if $\K$ and $\K'$ are finite,
then they have the same characteristic in view of Proposition
\ref{characteristic}. Now, the desired result is a consequence of
Theorem 31 in \cite{St68} if $\K$ and $\K'$ are finite and from
Theorem 8.16 of \cite{BT73} otherwise.
\end{proof}

The preceding theorem can be used to decompose automorphisms of a
given Kac-Moody group in a product of automorphisms of five
specific kinds, as mentioned in the introduction. For the precise
definitions of these specific automorphisms and further comments
on them, we refer the reader to Section 9 of \cite{CM03b}.
\begin{cor}\label{KM over arb fields cor}
Let $G$ be a Kac-Moody group over a field $\K$ of at least 4
elements and associated with a generalized Cartan matrix $A$ of
indecomposable type. Then, any automorphism of $G$ which preserve
bounded subgroups can be written as a product of an inner, a sign,
a diagonal, a graph and a field automorphism. Furthermore, if $G$
is ``simply connected in the weak sense'' (see \cite{Ti87}, Remark
3.7(c) p. 550) and if moreover, either $\mathrm{char}(\K)=0$ or
every off-diagonal entry of the generalized Cartan matrix $A$ is
prime to $\mathrm{char}(\K)$, then the term `graph automorphism'
may be replaced by `diagram automorphism' in the preceding
statement.
\end{cor}
The proof goes along the lines of the proof of Theorem 2.7 of
\cite{CM03b} and is omitted here.

\section{Kac-Moody groups over finite fields}

Corollary~\ref{cor 2} of the introduction is a consequence of the
following two results.

\begin{theo}\label{main over finite fields}
Let $\mathfrak{E}$ be the collection of all twin root data arising
from Kac-Moody groups over finite fields of cardinality at least
4. Then any element of $\mathfrak{E}$ is $\mathfrak{E}$-rigid.
\end{theo}
\begin{proof} Given any two elements $\mathcal{D}$ and $\mathcal{D}'$ of
$\mathfrak{E}$, it is clear from Corollary \ref{G_f <=> finite}
that every isomorphism of $G^\mathcal{D}$ to $G^\mathcal{D'}$
preserves bounded subgroups. The result is thus a consequence of
Theorem \ref{main over arbitrary fields}.
\end{proof}

\begin{cor}
Let $G$ be a Kac-Moody group over a finite field $\K$ of at least
4 elements and associated with a generalized Cartan matrix $A$ of
indecomposable type. Then, any automorphism of $G$ can be written
as a product of an inner, a sign, a diagonal, a graph and a field
automorphism. Furthermore, if $G$ is ``simply connected in the
weak sense'' (see \cite{Ti87}, Remark 3.7(c) p. 550) and if
moreover, every off-diagonal entry of the generalized Cartan
matrix $A$ is prime to $\mathrm{char}(\K)$, then the term `graph
automorphism' may be replaced by `diagram automorphism' in the
preceding statement.
\end{cor}
As for Corollary \ref{KM over arb fields cor} above, the proof is
omitted.

\end{document}